\documentclass[12pt]{amsart}
\usepackage{amsmath,amssymb,txfonts}
\usepackage{amssymb}
\usepackage{amsmath}
\usepackage{mathrsfs}
\usepackage{amsmath,amssymb}
\usepackage{amsmath}

\allowdisplaybreaks[4]
\usepackage[pagewise]{lineno}

\newtheorem{theorem}{Theorem}[section]
\newtheorem{lemma}[theorem]{Lemma}
\newtheorem{proposition}[theorem]{Proposition}

\theoremstyle{definition}
\newtheorem{definition}[theorem]{Definition}

\theoremstyle{remark}
\newtheorem{remark}[theorem]{Remark}

\def\R{\mbox{\rm{I}}\!\mbox{\rm{R}}}

\def\<{\leq}             \def\>{\geq}

\numberwithin{equation}{section}

\usepackage{color}

\begin{document}

\title[Naiver-Stokes equations]
{Well-posedness of Navier-Stokes equations established by the decaying speed of single norm}


\author[Qixiang Yang]{Qixiang Yang}
\address{School of Mathematics and Statistics, Wuhan University, Wuhan, 430072, China.}
\email{qxyang@whu.edu.cn}


\author[Huoxiong Wu]{Huoxiong Wu}
\address{School of Mathematical sciences, Xiamen University, Xiamen Fujian, 361005, China.
}
\email{huoxwu@xmu.edu.cn}

\author[Jianxun He]{Jianxun He}
\address{School of Mathematics and Information Sciences, Guangzhou University, Guangzhou, 510006, China.}
\email{hejianxun@gzhu.edu.cn}

\author[Zhenzhen Lou]{Zhenzhen Lou}
\address{School of Mathematics and Statistics, Qujing Normal University, Qujing, 655011, China \&
School of Mathematics and Information Sciences, Guangzhou University,
Guangzhou 510006, China}
\email{zhenzhenlou@e.gzhu.edu.cn}
\address{Corresponding author: Jianxun He, Zhenzhen Lou}

\thanks{Project supported by NSFC  Nos. 11571261, 12071229, 11771358, 11871101.}


\date{}

\dedicatory{}


\begin{abstract}
The decaying speed of a single norm
more truly reflects the intrinsic harmonic analysis structure of the solution of
the classical incompressible Navier-Stokes equations.
No previous work has been able to establish the well-posedness
under the decaying speed of a single norm with respect to time,
and the previous solution space is contained in the intersection of two spaces defined by different norms.
In this paper, for some separable initial space $X$,
we find some new solution space which is not the subspace of $L^{\infty}(X)$.
We use parametric Meyer wavelets to establish the well-posedness
via the decaying speed of a single norm only,
without integral norm to $t$.
\end{abstract}

\maketitle



 \vspace{0.1in}

\section{Introduction} 
\label{intro}

The Cauchy problem of Navier-Stokes equations on the half-space $\mathbb{R}^{1+n}_{+}=
(0,\infty) \times \mathbb{R}^{n}, n\geq 2,$ is defined as:
\begin{equation}\label{eqn:ns}
\left\{\begin{array}{ll} \partial_{ t} u
-\Delta u + u \cdot \nabla u -\nabla p=0,
& \mbox{ in } \mathbb{R}^{1+n}_{+}; \\
\nabla \cdot u=0,
& \mbox{ in } \mathbb{R}^{1+n}_{+}; \\
u|_{t=0}= a, & \mbox{ in } \mathbb{R}^{n}.
\end{array}
\right.
\end{equation}
Upon letting $R_{j}, j=1,2,\cdots n$ be the Riesz transforms,
writing
\begin{equation*}
\begin{cases}
\mathbb{P}= \{\delta_{l,l'}+ R_{l}R_{l'}\}, l,l'=1,\cdots,n;\\
\mathbb{P}\nabla (u\otimes u)= \sum\limits_{l}
\partial x_{l} (u_{l}u) -
\sum\limits_{l} \sum\limits_{l'} R_{l}R_{l'} \nabla (u_{l} u_{l'});\\
\widehat{e^{t\Delta}f}(\xi) =
e^{-t|\xi|^{2}}\hat{f}(\xi).
\end{cases}
\end{equation*}
And using $\nabla \cdot u=0$, we can see that solutions of the above
Cauchy problem are then obtained via the integral equation:
\begin{equation}\label{eqn:mildsolution}
\begin{cases}
u(t,x)= e^{t\Delta} a(x) - B(u,u)(t,x);\\
B(u,u)(t,x)\equiv\int^{t}_{0} e^{(t-s)\Delta}
\mathbb{P}\nabla (u\otimes u) ds,
\end{cases}
\end{equation}
which can be solved by a fixed-point method whenever the convergence
is suitably defined in some suitable function space. Solutions of
(\ref{eqn:mildsolution}) are called mild solutions of
(\ref{eqn:ns}). The notion of such a mild solution was pioneered by
Kato-Fujita \cite{KT} in 1960s.
For initial value spaces $X^{n}$, usually, one finds solution in $L^{\infty}(X^{n})$.
In this paper, {\bf our solution space defined by the decaying speed of a single norm is not a subspace of $L^{\infty}(X^{n})$}.

During the latest decades, many important results about mild solutions to (\ref{eqn:ns})
have been established.
A large amount  of work has established the well-posedness of
Navier-Stokes in different types of initial value spaces with two different norms.
See for example, Cannone \cite{C1, C2},
Germin-Pavlovic-Staffilani \cite{GPS}, Giga-Miyakawa \cite{GM},
Kato \cite{Kat},  Koch-Tataru \cite{KT},
Li-Xiao-Yang \cite{LXY}, Li-Yang \cite{YL},
Lin-Yang \cite{LY},
Miao-Yuan-Zhang \cite{C. Miao B. Yuan  B. Zhang 1},
Wu \cite{W1,W2,W3,W4}.
There are usually two situations:
\begin{itemize}
\item[\rm(1)]  For the first situation, the first norm is to take the norm of some index space variable related to the initial space and then take the infinite norm of time;
the second norm is to take the norm of another index space variable related to the initial space and then take the integral norm of time.
\item[\rm(2)] For the second situation, the first norm is the norm of some index space variable related to the initial space and then the infinite norm of time;
the second norm is the norm of another index space variable independent of the initial space and then the infinite norm of time.
\end{itemize}

Here we change the rule.
In this paper, we use parameter Meyer wavelets to introduce
some new spaces defined by the decaying speed of single norm.
For Besov space with negative regularity index,
by applying the new decomposition methods,
we can select the solution space
which has no relation to the integral norm of $t$
and which is defined only by the decaying speed of time $t$ of single norm.
{\bf After long-term calculation, we finally find that
the parameter wavelet structure can match the velocity of the norm decay of Besov space where $p$ is greater than the dimension $n$.}
For Besov spaces $\dot{B}^{\frac{n}{p}-1,p}_{p}$,
there exist some spaces $Y^{p}_{m}$ (see Definition 3.2) which are not subspaces of $L^{\infty}(\dot{B}^{\frac{n}{p}-1,p}_{p})$.
Based on this, we establish the following global well-posedness of
mild solutions for the Navier-Stokes equations with small
initial values in Besov  spaces $\dot{B}^{\frac{n}{p}-1,p}_{p}$.
\begin{theorem}\label{mthmain} Given $n<p<\infty$ and $ m\geq 1$.
The Navier-Stokes equations have a unique smooth solution in $(Y^{p}_{m})^n$ for all initial data $a(x)$ with $\nabla \cdot a =0$ and $\|a\|_{(\dot{B}^{\frac{n}{p}-1,p}_{p})^n}$ small enough.
\end{theorem}

\begin{remark}
(i) The purpose of this paper is to show the technique of establishing well-posedness of equations with the decaying speed of a single norm.
{\bf For simplicity}, we pick Besov space as the initial values.
Another advantage of this approach is that we can consider the well-posedness in some complex initial value spaces.

(ii) Triebel-Lizorkin-Lorentz spaces $\dot{L}^{s,q}_{p,\tau}$ and Besov-Lorentz spaces $\dot{B}^{s,q}_{p,\tau}$
were introduced in \cite{YCP}.
They have relation to the real interpolation spaces of Triebel-Lizorkin spaces and Besov spaces.
In fact, for $p>1, 0<\tau\leq \infty$, $\dot{L}^{0,2}_{p,\tau}$ are just the usual Lorentz spaces $L^{p,\tau}$.
The traditional methods need to consider the integration of the norm of the spatial variables with respect to time,
so it involves the problem of exchangeability between the integral with respect to time $t$ and the Lorentz index $\tau$.
 Therefore the traditional methods cannot deal with these spaces.
When Barraza \cite{Ba} considered the relative well-posedness of Lorentz spaces $L^{n,\infty}$,
he avoided this problem by limiting the initial values to functions that were homogeneous with degree -1
and he considered $L^{n}(S^{n-1})$.
Using the techniques in this paper,
we are studying the well-posedness of the initial values in general Lorentz spaces.
\end{remark}

\begin{remark}
Neededless to say, our current work grows from the already-known results.

(i) We have analyzed all existing results on the well-posednesss of the classical Navier-Stokes equations.
A common feature is that
all the above works handled equally the low and high frequencies,
or decompose the low and high frequencies equally.
See Cannone \cite{C1, C2},
Germin-Pavlovic-Staffilani \cite{GPS}, Giga-Miyakawa \cite{GM},
Kato \cite{Kat},  Lemari\'e \cite{Lem},
Li-Xiao-Yang \cite{LXY},
Lin-Yang \cite{LY},
Miao-Yuan-Zhang \cite{C. Miao B. Yuan  B. Zhang 1},
Wu \cite{W1,W2,W3,W4}, Li-Yang \cite{YL}, etc.
Improving the structure of solution space is one of our efforts
to improve the algorithm of mild solution to the extreme in the sense of harmonic analysis.
Our solution space $(Y^{p}_{m})^n$ is the first space defined with the decaying speed of a single norm.

(ii) For initial data spaces $X$, almost all the solution spaces are contained in $L^{\infty}(X)$.
Recently, Yang-Yang-Wu \cite{YYW} consider carefully Koch-Tataru's initial data space $({\rm BMO}^{-1})^{n}$ in \cite{KT}.
Yang-Yang-Wu found that Koch-Tataru's solution space $Y(({\rm BMO}^{-1})^{n})$
is not defined as a subspace of $L^{\infty}(({\rm BMO}^{-1})^{n})$.
Our solution space $(Y^{p}_{m})^n$ is the second space to have such property.
Except for the first difference mentioned above (i),
another diffence is that ${\rm BMO}^{-1}$ is an inseparable Banach  space,
but our initial data spaces are separable Banach spaces.

(iii) For the fractional Navier-Stokes equations, Yang-Yang \cite{YY} have noticed this phenomenon and
they found that there exists excessive decomposition for low frequency.
They have introduced parameter Meyer wavelets to avert the excessive decomposition of low frequency.
But there are some essential differences between the fractional Navier-Stokes equations
and classical Navier-Stokes  equations in the technique of dealing with the norm structure.
The former has well-posedness for Bloch space, see \cite{YY}.
For the later, until now, the biggest initial data space is Koch-Tataru's ${\rm BMO}^{-1}$,
and we did not know whether Bloch space has well-posedness or not.
For classical Navier-Stokes  equations, there is no solution space
defined by the decaying property of a single norm with respect to time
in all the previous works.
For $1\leq p\leq n$, we can not use  the new time-frequency structure in \cite{YY} to deal with the heat flow.
We have taken many efforts to apply parameter Meyer wavelet basis to introduce the present solution spaces.
\end{remark}

Our main ingredients are the following two folds:
(i). New time-space view based on parameter Meyer wavelets;
(ii). New solution spaces based on new harmoinc analysis structure defined by the decaying speed of a single norm only.
In the end of this section, we introduce the structure of the rest of the article.

In Section 2, we introduce some preliminary knowledge, namely the discretization technique and the basic inequalities.
Parameter Meyer wavelets no longer handle low and high frequencies equally as previous methods and
provide us a new time-space views for classic Navier-Stokes equations.
These new spatio-temporal methods of harmonic analysis will provide an effective tool for
the study of operators, classical functions and heterogeneous spaces later in this paper.

In Section 3, we use parameterized Meyer wavelets to analyze the heat flow and
we introduce adaptive evolution space $Y^{p}_{m}$ as the basis for constructing the solution space.
For Besov space $\dot{B}^{\frac{n}{p}-1,p}_{p}$ with negative regularity index,
by applying the new decomposition methods,
$Y^{p}_{m}$ has no relation to the integral norm of and
which is defined only by the single norm with respect to time.
In Remark 3.3, we mention that $Y^{p}_{m}$ space allows the unbounded norm of Besov space and
the discontinuity of Riesz operator on such a space,
and  we use the technique of acting Riesz operators on wavelets to prevent things like this from happening.

In Section 4, we select the solution space based on $Y^{p}_{m}$ and
we use the parameterized space-time view to treat the boundedness of the bilinear operators to the discrete case.
Unlike the traditional way of dealing with the classical Navier-Stokes equations,
we do not decompose the low frequencies.

In Section 5, we prove the continuity of discretized bilinear operators.
The main harmonic analysis techniques are
(1) to study the winding of Fourier transform of the wavelet functions,
(2) to consider the winding between the absolute value of the parameter wavelet coefficients and
the fast attenuation with respect to the discrete spatial variable and the discrete frequency,
(3) to estimate the corresponding  H\"{o}lder inequality.

With the help of this  new harmonic analysis solution structure,
we establish the global well-posedness of mild solutions for the Navier-Stokes equations with small initial values in Besov space.
In Section 6, first, we provide some comments on previous results and improvements on our current results.
Then we use the conclusions proved earlier in this paper to show the main theorem of this paper, namely well-posedness of the Navier-Stokes equations.

\section{Preliminaries}\label{sec2}

In this section, we introduce some preliminary knowledge,
namely the discretization technique and the basic inequalities.
These new spatio-temporal methods of harmonic analysis
provide an effective tool for the study of operators,
classical functions and heterogeneous spaces.

\subsection{Kernels and basic inequalities}
For all $t>0$ and $ l,l',l''\in \{1,\cdots, n\}$,
denote
$$A^{t}_{l}= e^{t\Delta} \partial x_{l},$$
$$A^{t}_{l,l',l''}=e^{t\Delta} \partial_{x_l}\partial x_{l'}\partial x_{l''}(-\Delta)^{-1}.$$
For $t=1$, let the relative kernels of these operators be $g_{l}$ and $g_{l,l',l''}$.

\begin{lemma}\label{kernel} $g_{l}$ and $g_{l,l',l''}$ satisfy
\begin{equation}\label{el}(1+|x|)^{n+1} |g_{l}(x)|\leq C,\end{equation}
\begin{equation}\label{ell}(1+|x|)^{n+1} |g_{l,l',l''}(x)|\leq C.\end{equation}
\end{lemma}

For any $ l,l',l''\in \{1,\cdots, n\}$,
denote
$$A^{t-s }_{l}  f
= \int (t-s)^{-\frac{n+1}{2}} g_{l}\left(\frac{x-y}{(t-s)^{\frac{1}{2}}}\right) f(y) dy,$$
$$A^{t-s }_{l,l',l''} f
= \int (t-s)^{-\frac{n+1}{2}} g_{l,l',l''}\left(\frac{x-y}{(t-s)^{\frac{1}{2}}}\right) f(y) dy.$$

The following inequality can be found in \cite{Me}:
\begin{lemma}\label{le:2.2}
If $a\geq 1$, then
$$1+|x| \leq 2(1+|x-y|)(1+a|y|).$$

\end{lemma}

\begin{proof}
We distinguish two cases: (i) $a|y|\leq \frac{|x|}{2}$ and (ii) $a|y|\geq \frac{|x|}{2}$.
We get the above conclusion.
\end{proof}

By the above Lemma \ref{le:2.2}, we have:
\begin{lemma}\label{le:2.3}
For $ j\geq j'$, we have
$$\begin{array}{rl}
&(1+|2^{j'}x-k'|)^{-n-1} (1+ |2^{j}x-k|)^{-N-n-1} \\
\lesssim &   (1+ |2^{j'-j}k-k'|)^{-n-1} (1+ |2^{j} x-k|)^{-N}.
\end{array}$$
\end{lemma}


\subsection{Meyer wavelets} First of all,
we indicate that we will use tensorial product real valued
orthogonal Meyer wavelets. We refer the reader to \cite{Me}, \cite{Woj} and
\cite{Yang1} for further information.
Let $\Psi^{0}$ be an even function in $ C^{\infty}_{0}
([-\frac{4\pi}{3}, \frac{4\pi}{3}])$ with
\begin{equation}
\left\{ \begin{aligned}
&0\leq\Psi^{0}(\xi)\leq 1; \nonumber\\
&\Psi^{0}(\xi)=1\text{ for }|\xi|\leq \frac{2\pi}{3}.\nonumber
\end{aligned} \right.
\end{equation}
  Write
  $$\Omega(\xi)= \sqrt{(\Psi^{0}(\frac{\xi}{2}))^{2}-(\Psi^{0}(\xi))^{2}}.$$ Then $\Omega(\xi)$ is an even
function in $ C^{\infty}_{0}([-\frac{8\pi}{3}, \frac{8\pi}{3}])$.
Clearly,

\begin{equation}
\left\{ \begin{aligned}
&\Omega(\xi)=0\text{ for }|\xi|\leq \frac{2\pi}{3};\nonumber\\
&\Omega^{2}(\xi)+\Omega^{2}(2\xi)=1=\Omega^{2}(\xi)+\Omega^{2}(2\pi-\xi)\text{
for }\xi\in [\frac{2\pi}{3},\frac{4\pi}{3}].\nonumber
\end{aligned} \right.
\end{equation}

 Let $\Psi^{1}(\xi)=
\Omega(\xi) e^{-\frac{i\xi}{2}}$. For any $\epsilon=
(\epsilon_{1},\cdots, \epsilon_{n}) \in \{0,1\}^{n}$, define
$\Phi^{\epsilon}(x)$ by $\hat{\Phi}^{\epsilon}(\xi)=
\prod\limits^{n}_{i=1} \Psi^{\epsilon_{i}}(\xi_{i})$. For $j\in
\mathbb{Z}$ and $k\in\mathbb{Z}^{n}$, let $\Phi^{\epsilon}_{j,k}(x)=
2^{\frac{nj}{2}} \Phi^{\epsilon} (2^{j}x-k)$. For all $\epsilon\in \{0,1\}^{n}, j\in \mathbb{Z}, k\in \mathbb{Z}^{n}$ and distribution $f(x)$, denote
$f^{\epsilon}_{j,k}=\langle f, \Phi^{\epsilon}_{j,k}\rangle$.
Furthermore, we put
\begin{equation*}
\Gamma =\{(\epsilon,k), \epsilon\in \{0,1\}^{n}\backslash\{0\}, k\in \mathbb{Z}^{n}\},\,
\Lambda =\{(\epsilon,j,k), (\epsilon,k)\in \Gamma, j\in\mathbb{Z}\}.
\end{equation*}

The following result is well-known.
\begin{lemma}\label{le1}
(i) The Meyer wavelets $\{\Phi^{\epsilon}_{j,k}(x)\}_{(\epsilon,j,k)\in
\Lambda}$ form an  orthogonal basis in $L^{2}(\mathbb{R}^{n})$.
For any $f\in L^{2}(\mathbb{R}^{n})$, the following wavelet
decomposition holds in the $L^2$ convergence sense:
$$\begin{array}{c}
f(x)=\sum\limits_{(\epsilon,j,k)\in\Lambda}f^{\epsilon}_{j,k}\Phi^{\epsilon}_{j,k}(x).
\end{array}$$

(ii) Given $s\in \mathbb{R}$ and $1\leq p,q\leq \infty$.
$f(x)\in \dot{B}^{s,q}_{p} \Longleftrightarrow$
$$\begin{array}{c}
\sum\limits_{j}  2^{jq(s+\frac{n}{2}-\frac{n}{p})}  \left( \sum\limits_{(\epsilon, k)\in\Lambda}|f^{\epsilon}_{j,k}|^{p}\right)
^{\frac{q}{p}}<\infty.
\end{array}$$
\end{lemma}

\subsection{Parameter wavelets and relative inequalities}
To adapt the heat flow, we use parameter Meyer wavelets. That is to say, for different $t$, we use different wavelet basis
and we will prove these wavelets can be adapted to the study of heat flow  in  Section \ref{sec4}.
For $t>0$, denote $j_t$ the smallest integer such that $2^{2j} t\geq 1$.
If $t$ satisfies $1\leq 2^{2 j_t} t< 4$, then we denote $t\in S_{j_t}$.
Denote $\Lambda_{t}^{(0)}= \{ (0,j_t,k), k\in \mathbb{Z}^n\}$,
$\Gamma = \{ (\epsilon,k), \epsilon=\{0,1\}^n \backslash \{0\},  k\in \mathbb{Z}^n\}$ and
$\Lambda^{(1)}_{t} = \{ (\epsilon,j,k), \epsilon=\{0,1\}^n \backslash \{0\}, j\in \mathbb{Z}, k\in \mathbb{Z}^n, j\geq j_t\}$.
Let $\Lambda_{t}= \Lambda^{(0)}_{t} \bigcup \Lambda^{(1)}_{t}$.

\begin{lemma}\label{le:2.5} For $t>0$,
$\{\Phi^{\epsilon}_{j,k}\}_{(\epsilon,j,k)\in \Lambda_{t}}$ are orthogonal wavelet basis in $L^{2}(\mathbb{R}^{n})$.
\end{lemma}

When relative operators act on the parameter wavelets, by Lemmas \ref{kernel} and \ref{le:2.3}, we have the following estimations:
\begin{lemma}
Given $\alpha\in \mathbb{N}^{n}, N\geq n+1, k\in \mathbb{Z}^{n}$ and $0\leq s\leq t$.
$$\begin{array}{c} \left|e^{(t-s)\Delta}\partial^{\alpha} \Phi^{0}_{j_t,k}\right| \lesssim  2^{(\frac{n}{2}+|\alpha|)j_t} \left(1+|2^{j_t}x-k|\right)^{-N},\, \forall\, 0\leq |\alpha| \leq 2.\\
\left|e^{(t-s)\Delta}\partial^{\alpha} (-\Delta)^{-1} \Phi^{0}_{j_t,k}\right| \lesssim  2^{(\frac{n}{2}+|\alpha|)j_t} \left(1+|2^{j_t}x-k|\right)^{-n-|\alpha|+2 }, \;\forall\; 3\leq |\alpha|\leq 5.
\end{array}$$

For $j\geq j_t, \epsilon\neq 0$ and $ N\geq n+1$, there exists $c>0$ such that
$$ \begin{array}{c} \left|e^{(t-s)\Delta}\partial^{\alpha} \Phi^{\epsilon}_{j,k}\right| \lesssim  e^{-c(t-s) 2^{2j}} 2^{(\frac{n}{2}+|\alpha|)j} \left(1+|2^{j}x-k|\right)^{-N},\; \forall\; 0\leq |\alpha|\leq 2,\\
\left|e^{(t-s)\Delta}\partial^{\alpha} (-\Delta)^{-1} \Phi^{\epsilon}_{j,k}\right| \lesssim  e^{-c(t-s) 2^{2j}} 2^{(\frac{n}{2}+|\alpha|)j} \left(1+|2^{j}x-k|\right)^{-N}\!\!\!,\!\!\;\forall \;3\leq |\alpha|\leq 5.
\end{array}$$
\end{lemma}

Further, for all $\epsilon,\epsilon'\in \{0,1\}^{n}\backslash \{0\}$ and $k,k',k''\in \mathbb{Z}^{n}$, denote:
\begin{equation}\label{ea}
a^{\epsilon,\epsilon'}_{j_t,k,k',k''}(t)=
\left\langle \sum\limits_{j\geq 2+j_s}  \Phi^{\epsilon}_{j,k}(x) \Phi^{\epsilon'}_{j,k''}(x),
A^{t-s}_{l,l',l''} \Phi^{0}_{j_t, k'}\right\rangle.
\end{equation}
Applying the above inequalities (\ref{el}), (\ref{ell}) in Lemma \ref{kernel}
and applying Lemma \ref{le:2.3}, we have
\begin{lemma}\label{lem:es1} $a^{\epsilon,\epsilon'}_{j_t,k,k',k''}(t)$ in \eqref{ea} satisfies
$$\left|a^{\epsilon,\epsilon'}_{j_t,k,k',k''}(t)\right| \lesssim  2^{j_t+\frac{nj_t}{2}} \sum\limits_{j\geq 2+j_s}   \left(1+|k-k''|\right)^{-N} \left(1+|2^{j_t-j}k-k'|\right)^{-n-1}.$$
\end{lemma}

\section{Adaptive evolution space}\label{sec4}
The Navier-Stokes equations
(\ref{eqn:ns}) are invariant under the scaling
\begin{equation*}
\begin{cases}
u_{\lambda}(t,x) = \lambda u(\lambda^{2}t, \lambda x);\\
p_{\lambda}(t,x) = \lambda^{2} p(\lambda^{2}t, \lambda x),
\end{cases}
\end{equation*}
If $u(t,x)$ is a solution of (1.1) and we replace
$u(t,x), p(t,x), a(x)$ by $$u_{\lambda}(t,x) = \lambda
u(\lambda^{2}t, \lambda x),\ p_{\lambda}(t,x) = \lambda^{2} p(\lambda^{2}t, \lambda x)$$ and
$a_{\lambda}(x) = \lambda a(\lambda x),$ respectively, then
$u_{\lambda}(t,x)$ is also a solution for (1.1).
If a spaces $X$ satisfying that $\|a\|_{X}\sim \|a_{\lambda}(x)\|_{X}$,
then $X$ is said to be a critical space.
Critical spaces occupied a significant place for Navier-Stokes equations (\ref{eqn:ns}).

Usually, the solution space $Y(X)$ is chosen to be the subspace of $L^{\infty}(X)$.
And all the previous solution spaces $Y(X)$ have not been defined by the decay speed of a single norm.
See Cannone \cite{C1, C2},
Germin-Pavlovic-Staffilani \cite{GPS}, Giga-Miyakawa \cite{GM},
Kato \cite{Kat},  Li-Xiao-Yang \cite{LXY}, Lin-Yang \cite{LY},
Miao-Yuan-Zhang \cite{C. Miao B. Yuan  B. Zhang 1},
Wu \cite{W1,W2,W3,W4} and Yang-Li \cite{YL}.
To simplify the notations, we consider the critical Besov space $\dot{B}^{\frac{n}{p}-1,p}_{p}$,
but our solution space is not a subspace of $L^{\infty}((\dot{B}^{\frac{n}{p}-1,p}_{p})^{n}).$
In this section, we prove first that parameter Meyer wavelets are adapted to the heat flow.
Then we apply such idea to establish the relation between Besov space $\dot{B}^{\frac{n}{p}-1,p}_{p}$
and their relative adaptive evolution space $Y^{p}_{m}$.

\subsection{Heat flow}

Assume that $f(x)\in S'(\mathbb{R}^{n})$ and $f(t,\cdot)\in S'(\mathbb{R}^{n})$.
For Meyer wavelets
$\{\Phi^{\epsilon}_{j,k}\}_{(\epsilon,j,k)\in\Lambda}$, let
$a^{\epsilon}_{j,k} =
\langle f, \Phi^{\epsilon}_{j,k}\rangle$.
For parameter Meyer wavelets
$\{\Phi^{\epsilon}_{j,k}\}_{(\epsilon,j,k)\in\Lambda_t}$, let
$a^{\epsilon}_{j,k}(t) = \langle f(t,x),
\Phi^{\epsilon}_{j,k}\rangle$. By Lemmas \ref{le1} and \ref{le:2.5},
$$\begin{array}{rl}
&f(x)= \sum\limits_{(\epsilon,j,k)\in \Lambda} a^{\epsilon}_{j,k}
\Phi^{\epsilon}_{j,k}(x) \text{ and }  f(t,x)=
\sum\limits_{(\epsilon, j, k)\in \Lambda_{t}} a^{\epsilon}_{j,k}(t)
\Phi^{\epsilon}_{j,k}(x).
\end{array}$$

If $f(t,x)= e^{t\Delta}f(x)$, then $\forall (\epsilon,j,k)\in \Lambda^{(0)}_{t}$, we have
\begin{equation}\label{eq3}
\begin{array}{rl}
a^{0}_{j_t,k}(t)
= &\sum\limits_{\epsilon',j'\leq 1+j_t, k'}
a^{\epsilon'}_{j',k'}\left \langle e^{t\Delta}
\Phi^{\epsilon'}_{j',k'}(x),
\Phi^{0}_{j_t,k}(x)\right\rangle\\
= &\sum\limits_{\epsilon',j'\leq 1+j_t, k'}
a^{\epsilon'}_{j',k'}\left \langle
\Phi^{\epsilon'}_{j',k'}(x),
e^{t\Delta}\Phi^{0}_{j_t,k}(x)\right\rangle.
\end{array}
\end{equation}
$For\;all\; (\epsilon,j,k)\in \Lambda^{(1)}_{t}$, we obtain
\begin{equation}\label{eq4}
\begin{array}{rl}
\,\,\,a^{\epsilon}_{j,k}(t)
=\!\! &\sum\limits_{\epsilon',|j-j'|\leq 1, k'}
a^{\epsilon'}_{j',k'} \left\langle e^{t\Delta}
\Phi^{\epsilon'}_{j',k'}(x),
\Phi^{\epsilon}_{j,k}(x)\right\rangle\\
=\!\! &\sum\limits_{\epsilon',|j-j'|\leq 1, k'}
a^{\epsilon'}_{j',k'} \left\langle
\Phi^{\epsilon'}_{j',k'}(x),
e^{t\Delta} \Phi^{\epsilon}_{j,k}(x)\right\rangle.
\end{array}
\end{equation}

By (\ref{eq3}) and (\ref{eq4}), we can see that parameter Meyer wavelets are adapted to heat flow.
\begin{lemma}\label{le4}
There exist a  constant $N>n$ large
enough and a fixed small constant $ c>0$ such that
\begin{itemize}
\item[\rm(i)]
if $ (\epsilon,j,k)\in \Lambda^{(1)}_{t}$, then
\begin{equation}\label{sw}
\begin{array}{l}
\left|a^{\epsilon}_{j,k}(t)\right|\lesssim e^{- c t 2^{2j}}
\sum\limits_{\epsilon',|j-j'|\leq 1, k'} \left|a^{\epsilon'}_{j',k'}\right|
\left(1+\left|2^{j-j'}k'-k\right|\right)^{-N};
\end{array} \end{equation}
\item[\rm(ii)]if $ (\epsilon,j,k)\in \Lambda^{(0)}_{t}$ , then
\begin{equation}\label{eq:3.4}
\begin{array}{rl}\left|a^{0}_{j_t,k}(t)\right|\lesssim  \sum\limits_{j'\leq 1+j_t} 2^{\frac{n(j'-j_t)}{2}}
\sum\limits_{\epsilon',k'} \left|a^{\epsilon'}_{j',k'}\right|
\left(1+\left|k'-2^{j'-j_t}k\right|\right)^{-N}. \end{array}\end{equation}
\end{itemize}
\end{lemma}

\subsection{Adaptive evolution space}\label{Y}
Traditionally, the norm of  solution space consists of two parts.
The solution is in the intersection of two kinds of spaces.
In this paper, we introduce new solution structures.
Applying the above adaptive wavelets, we introduce adaptive evolution space
$Y^{p}_{m}$ related to the Besov space $\dot{B}^{\frac{n}{p}-1,p}_{p}$.

The method we adopt is to divide first time-space into binary ring time.
Then in the ring time, the low frequency matching with the ring time is not decomposed,
and the high frequency is decomposed.
For each ring frequency, the Lebesgue norm for the spatial variable is first obtained,
and then the upper bound on the ring time is obtained.
These discrete quantities are first summed in terms of location and frequency, and then bounded in discrete time.

\begin{definition} Given $0<m,p<\infty.$ $u(t,x)\in Y^{p}_{m}$ if and only if
$$\sup\limits_{j_t\in \mathbb{Z}}\sum\limits_{j\geq j_t}2^{2mp(j-j_t)} 2^{jp(\frac{n}{2}-1)}
\sup\limits_{1\leq t 2^{2j_t}<4} \sum\limits_{\epsilon, k:(\epsilon,j,k)\in \Lambda_t}\left|a^{\epsilon}_{j,k}(t)\right|^{p}<\infty.$$
\end{definition}

For convenience, we introduce the following notations.

For all $j',j\in \mathbb{Z}, j\geq j'$ and $1\leq t 2^{2j'}<4 $, denote
$$\begin{array}{c}
A^{p}_{j'}=  2^{pj'(\frac{n}{2}-1)} \sup\limits_{1\leq t 2^{2j'}<4} \sum\limits_{k\in \mathbb{Z}^{n}} \left|a^{0}_{j',k}(t)\right|^{p} ,\\
A^{p,m}_{j,j'}= 2^{2mp(j-j')} 2^{pj(\frac{n}{2}-1)} \sup\limits_{1\leq t 2^{2j'}<4} \sum\limits_{(\epsilon,k)\in \Gamma}\left|a^{\epsilon}_{j,k}(t)\right|^{p},
\end{array}$$
$$\begin{array}{c}
H^{p}_{0}= \sup\limits_{j_t\in \mathbb{Z}}\left \{\sup\limits_{1\leq t 2^{2j_t}<4} 2^{pj_t(\frac{n}{2}-1)} \sum\limits_{k\in \mathbb{Z}^{n}} \left|a^{0}_{j_t,k}(t)\right|^{p}\right\}^{\frac{1}{p}}= \left\{\sup\limits_{j'} A^{p}_{j'}\right\}^{\frac{1}{p}},\\
H^{p}_{m}=\sup\limits_{j_t\in \mathbb{Z}}\left\{\sum\limits_{j\geq j_t}2^{2mp(j-j_t)} 2^{jp(\frac{n}{2}-1)} \sup\limits_{1\leq t 2^{2j_t}<4} \sum\limits_{(\epsilon,k)\in \Gamma}\left|a^{\epsilon}_{j,k}(t)\right|^{p}  \right\}^{\frac{1}{p}}
=\sup\limits_{j'} \left\{\sum\limits_{j\geq j'} A^{p,m}_{j,j'}\right\}^{\frac{1}{p}}\!\!\!.
\end{array}$$
Hence $u(t,x)\in Y^{p}_{m}$ if and only if $ H^{p}_{0} <\infty$ and $H^{p}_{m}<\infty.$

\begin{remark}
(i) For Besov spaces $\dot{B}^{\frac{n}{p}-1,p}_{p}$,
the solution spaces studied in \cite{C2, Lem, LXY, LY, W3} and \cite{YL}
are contained in $L^{\infty}((\dot{B}^{\frac{n}{p}-1,p}_{p})^{n})$.
Adaptive evolution space $Y^{p}_{m}$ is different to the above solution spaces.
$Y^{p}_{m}$ is not contained in $L^{\infty}(\dot{B}^{\frac{n}{p}-1,p}_{p})$ and
$Y^{p}_{m}$ more accurately reflects the internal harmonic analysis structure of the solution itself.

(ii) Riesz operators are not bounded on parameter Besov spaces $Y^{p}_{m}$.
Instead of studying the continuity of Riesz operators,
we directly consider the effect of Riesz operators on wavelets.
\end{remark}

We can prove that heat flow maps Besov spaces into
the above adaptive evolution space.
\begin{theorem} \label{th:B-to-Y}
Given  $n<p<\infty$ and $m>0 $.
If $f\in \dot{B}^{\frac{n}{p}-1,p}_{p}$, then $ e^{t\Delta} f \in Y^{p}_{m}.$
\end{theorem}

\begin{proof}
We apply the wavelet characterization of $\dot{B}^{\frac{n}{p}-1,p}_{p}$ in Lemma \ref{le1} to prove this Theorem.

(I) Firstly, we consider  $(0,j_t,k)\in \Lambda^{(0)}_{t}$. By \eqref{eq:3.4},
$$\begin{array}{rcl}
\left|a^{0}_{j_t,k}(t)\right| &\lesssim &  \sum\limits_{j'\leq -1+j_t} 2^{\frac{n(j'-j_t)}{2}}
\sum\limits_{\epsilon',k'} \left|a^{\epsilon'}_{j',k'}\right| \left(1+\left|k'-2^{j'-j_t}k\right|\right)^{-2N}\\
&\lesssim &  \sum\limits_{j'\leq -1+j_t} 2^{\frac{n(j'-j_t)}{2}}
\left[\sum\limits_{\epsilon',k'} \left|a^{\epsilon'}_{j',k'}\right|^{p} \left(1+\left|k'-2^{j'-j_t}k\right|\right)^{-2N}\right]^{\frac{1}{p}}.
\end{array}$$
Take small positive $\delta$ such that $0<\delta <\frac{p-n}{p}$, we have
$$\begin{array}{rcl}
\left|a^{0}_{j_t,k}(t)\right|
&\lesssim &  \left[\sum\limits_{j'\leq -1+j_t} 2^{ p (j'-j_t) (\frac{n}{2}-\delta) }
\sum\limits_{\epsilon',k'} \left|a^{\epsilon'}_{j',k'}\right|^{p} \left(1+\left|k'-2^{j'-j_t}k\right|\right)^{-2N}\right]^{\frac{1}{p}}.
\end{array}$$
Hence
$$\begin{array}{rcl}
H^{p}_{0} &\lesssim & 2^{j_t( \frac{n}{2}-1)} \left[\sum\limits_{k} \sum\limits_{j'\leq -1+j_t} 2^{ p (j'-j_t) (\frac{n}{2}-\delta) }
\sum\limits_{\epsilon',k'} \left|a^{\epsilon'}_{j',k'}\right|^{p} \left(1+\left|k'-2^{j'-j_t}k\right|\right)^{-2N}\right]^{\frac{1}{p}}\\
&\lesssim &   2^{j_t( \frac{n}{2}-1)} \left(\sum\limits_{j'\leq -1+j_t} 2^{ p (j'-j_t) (\frac{n}{2}-\delta) } 2^{  -n(j'-j_t)  }
\sum\limits_{\epsilon',k'} \left|a^{\epsilon'}_{j',k'}\right|^{p} \right)^{\frac{1}{p}}   \\
&\lesssim &    \left(\sum\limits_{j'\leq -1+j_t} 2^{  (j'-j_t) (p(1-\delta)-n) }
2^{pj'( \frac{n}{2}-1)}\sum\limits_{\epsilon',k'} \left|a^{\epsilon'}_{j',k'}\right|^{p} \right)^{\frac{1}{p}}.
\end{array}$$
Since $p(1-\delta)>n$, we obtain
$$H^{p}_{0} \lesssim  \sum\limits_{j'\leq -1+ j_t } 2^{  (j'-j_t) (p(1-\delta)-n) }
\lesssim C.
$$

(II) Secondly, we  think about  $(\epsilon,j,k)\in \Lambda^{(1)}_{t}$ and $f\in \dot{B}^{\frac{n}{p}-1,p}_{p}$.
By applying (\ref{sw}), we get
\begin{equation*}
\begin{array}{rcl}
\left|a^{\epsilon}_{j,k}(t)\right| &\lesssim &e^{- c t 2^{2j}}
\sum\limits_{\epsilon',|j-j'|\leq 1, k'} \left|a^{\epsilon'}_{j',k'}\right|
\left(1+\left|2^{j-j'}k'-k\right|\right)^{-N}\\
&\lesssim &e^{- c t 2^{2j}}
\left[\sum\limits_{\epsilon',|j-j'|\leq 1, k'} \left|a^{\epsilon'}_{j',k'}\right|^{p}
\left(1+\left|2^{j-j'}k'-k\right|\right)^{-N}\right]^{\frac{1}{p}}.
\end{array} \end{equation*}
Hence
$$\begin{array}{rcl}
H^{p}_{m} &\lesssim &
\sum\limits_{j\geq j_t} 2^{2mp(j-j_t)} 2^{pj(\frac{n}{2}-1)}
e^{-p c t 2^{2j}}\\
&&\times \sum\limits_{(\epsilon,k)\in \Gamma}
\sum\limits_{\epsilon',|j-j'|\leq 1, k'} \left|a^{\epsilon'}_{j',k'}\right|^{p}
\left(1+\left|2^{j-j'}k'-k\right|\right)^{-N}\\
&\lesssim &
\sum\limits_{j\geq j_t} 2^{2mp(j-j_t)} 2^{pj(\frac{n}{2}-1)}
e^{-p c t 2^{2j}} \sum\limits_{\epsilon',|j-j'|\leq 1, k'} \left|a^{\epsilon'}_{j',k'}\right|^{p}\\
&\lesssim &
\sum\limits_{j\geq j_t} 2^{2mp(j-j_t)}
e^{-p c t 2^{2j}} \lesssim C.
\end{array}$$

\end{proof}

By definition of $Y^{p}_{m}$ and wavelet knowledge, we can easily get  the following property.
\begin{proposition}
Given $m>0$.
$f\in Y^{p}_{m}$ implies
$$\sup\limits_{k\in \mathbb{Z}^n } 2^{\frac{nj}{2}} \left|a^{0}_{j,k}(t)\right|
\lesssim  2^{j_t},\; for \;all\; j\geq j_t.$$
\end{proposition}

\section{Boundedness of bilinear operators}

In this section, we use the parameterized space-time view to transform the boundedness of the following bilinear operators.
For $l,l',l''=1,\cdots,n$, define
\begin{equation}\label{41eq}
B_{l}(u,v)= \int^{t}_{0} A^{t-s}_{l}(u(x,s) v(x,s)) ds,
\end{equation}
\begin{equation}\label{42eq}
B_{l,l',l''}(u,v)= \int^{t}_{0} A^{t-s}_{l,l',l''} (u(x,s) v(x,s)) ds.
\end{equation}

The operators defined in the above equations \eqref{41eq} and \eqref{42eq} are bounded:
\begin{theorem} \label{th:4.1} Given $n<p<\infty$ and $m\geq 1$. For $l,l',l''=1,\cdots,n$, we have
$$\begin{array}{c}
(u,v)\mapsto B_{l}(u,v) \mbox { are bounded from } Y^{p}_{m }\times Y^{p}_{m } \mbox{ to } Y^{p}_{m },\\
(u,v)\mapsto B_{l,l',l''}(u,v) \mbox { are bounded from } Y^{p}_{m }\times Y^{p}_{m } \mbox{ to } Y^{p}_{m }.
\end{array}$$
\end{theorem}

\begin{proof}
For $j\in\mathbb{Z}$ and $\epsilon\in \{0,1\}^{n}\backslash \{0\}$, let
$$\begin{array}{rcl}
P_{j}f(x)= \sum\limits_{k\in \mathbb{Z}^{n}} f^{0}_{j,k}
\Phi^{0}_{j,k}(x),&& \\
Q_{j}^{\epsilon}f(x)= \sum\limits_{k \in \mathbb{Z}^n
} f^{\epsilon}_{j,k} \Phi^{\epsilon}_{j,k}(x)
&\text{ and }& Q_{j}f(x)= \sum\limits_{\epsilon \in
\{0,1\}^n\backslash \{0\}} Q^{\epsilon}_{j} f(x). \end{array}$$
By using the Meyer wavelets in Section 2,  the product of any two functions $u$
and $v$ can be decomposed as
\begin{equation}\label{eq:decompose}
\begin{array}{rl}
u(t,x)v(t,x)=\!\! & \!\!\sum\limits_{j\geq 2+j_t} \big\{ P_{j-2}u Q_{j}v + Q_{j-2}u Q_j v + Q_{j-1}u Q_j v+ Q_j u Q_j v \\
&+ Q_j u Q_{j-1}v + Q_j u Q_{j-2}v + Q_j u P_{j-2}v\big\}\\
&+ P_{1+j_t}u Q_{1+j_t}v+ Q_{1+j_t}u Q_{1+j_t}v +Q_{1+j_t}u P_{1+j_t}v\\
&+ P_{j_t}u Q_{j_t}v+ Q_{j_t}u Q_{j_t}v +Q_{j_t}u P_{j_t}v + P_{j_t}u P_{j_t}v.
\end{array}
\end{equation}

For $\epsilon\in \{0,1\}^n\backslash \{0\}$,
denote
\begin{equation*}
\begin{array}{lcl}
I^{\epsilon}_{u,v}& = & \sum\limits_{j\geq 2+j_t}  P_{j-2}u Q^{\epsilon}_{j}v ,  \\
II_{u,v}&= & \sum\limits_{j\geq 2+j_t} Q_j u Q_{j}v, \\
III_{u,v}&=& \sum\limits_{k, k'} u^{0}_{j_t, k} v^{0}_{j_t, k'} \Phi^{0}_{j_t, k} \Phi^{0}_{j_t, k'}.
\end{array}
\end{equation*}

When considering the terms in the decomposition equality (\ref{eq:decompose}), there are some similarities:
\item{\rm(i)}
For $\sum\limits_{j\geq 2+j_t}  P_{j-2}u Q^{\epsilon}_{j}v$ and $\sum\limits_{j\geq 2+j_t}  P_{j-2}v Q^{\epsilon}_{j}u$, by similarity, we consider
$I^{\epsilon}_{u,v}$.

\item{\rm(ii)}
For $\sum\limits_{j\geq 2+j_t} Q_{j-2}u Q_j v $, $\sum\limits_{j\geq 2+j_t} Q_{j-1}u Q_j v$, $\sum\limits_{j\geq 2+j_t} Q_j u Q_j v $,
$\sum\limits_{j\geq 2+j_t} Q_j u Q_{j-1}v$ and $\sum\limits_{j\geq 2+j_t} Q_j u Q_{j-2}v$, by similarity, we consider only  $II_{u,v}$.

\item{\rm(iii)}
For the terms $P_{1+j_t}u Q_{1+j_t}v$, $ Q_{1+j_t}u Q_{1+j_t}v$, $Q_{1+j_t}u P_{1+j_t}v$,
$ P_{j_t}u Q_{j_t}v$, $ Q_{j_t}u Q_{j_t}v $, $ Q_{j_t}u P_{j_t}v$ and $ P_{j_t}u P_{j_t}v$, by similarity, we   consider only the last term.
That is to say, we consider $III_{u,v}$.

Further, $B_{l}(u,v)$ is more easy to prove. We consider only $B_{l,l',l''}(u,v)$.
By similarity of skills, we need only to prove the boundedness of the following three bilinear operators:
\begin{align*}
B_{l,l',l'',\epsilon}(u,v)&= \int^{t}_{0} A^{t-s}_{l,l',l''} \left(\sum\limits_{j\geq 2+j_s}  P_{j-2}u Q^{\epsilon}_{j}v\right) ds,\\
B_{l,l',l'',1}(u,v)&=\int^{t}_{0} A^{t-s}_{l,l',l''}\left(\sum\limits_{j\geq 2+j_s} Q_j u Q_{j}v\right) ds,\\
B_{l,l',l'',2}(u,v)&=\int^{t}_{0} A^{t-s}_{l,l',l''}\left(III_{u,v}\right) ds.
\end{align*}

That is to say, to prove Theorem \ref{th:4.1}, we prove the following lemma:
\begin{lemma} \label{Le:4} Given $m\geq 1$, $n<p<\infty$, $l,l',l''=1,\cdots,n\;and\; \epsilon\in \{0,1\}^{n}\backslash \{0\}$.
If $u,v\in Y^{p}_{m}$, then
\begin{align*}
\begin{array}{ll}
{\rm (i)}\,\, & B_{l,l',l'',\epsilon}(u,v) \in Y^{p}_{m },\\
{\rm (ii)}\,\, & B_{l,l',l'',1}(u,v) \in Y^{p}_{m },\\
{\rm (iii)}\,\, & B_{l,l',l'',2}(u,v) \in Y^{p}_{m }.
\end{array}
\end{align*}
\end{lemma}
We will prove this lemma in the following section.
\end{proof}

\section{ The proof  of Lemma \ref{Le:4}}
The proof of  Lemma \ref{Le:4} is divided into three subsections.
The main harmonic analysis techniques are (1) to study the winding of Fourier transform of the wavelet functions,
(2) to consider the winding between the absolute value of the parameter wavelet coefficients
with respect to the discrete spatial variable and the fast attenuation corresponding to the discrete frequency,
(3)  to estimate the corresponding  H\"{o}lder inequality.

\subsection{ Proof of (i) of Lemma \ref{Le:4}}
We consider first the case where $\epsilon=0$.
Note that
$${\rm Supp} \widehat{P_{j-2}u} \subset \Big\{\left|\xi_i \right|\leq \frac{\pi}{3} \cdot 2^{j}, \forall i=1,\cdots, n\Big\},$$
$$ {\rm Supp} \widehat{Q^{\epsilon}_{j} u} \subset
\Big\{|\xi_i|\leq \frac{4\pi}{3} \cdot 2^{j}, \mbox{ if } \epsilon_i=0;
\frac{2\pi}{3} \cdot 2^{j}\leq
|\xi_i|\leq \frac{8\pi}{3} \cdot 2^{j}, \mbox{ if } \epsilon_i=1\Big\}.$$
Hence the support of the Fourier transform of $P_{j-2}u Q^{\epsilon}_{j} v$ is contained in a ring:
$$\Big\{|\xi_i|\leq \frac{5\pi}{3} \cdot 2^{j}, \mbox{ if } \epsilon_i=0;
\frac{\pi}{3} \cdot 2^{j}\leq
|\xi_i|\leq 3\pi \cdot 2^{j}, \mbox{ if } \epsilon_i=1\Big\}.$$
This implies
$${\rm Supp} \widehat{I^{\epsilon}_{u,v}} \subset \Big\{|\xi_i|\geq \frac{4\pi}{3} \cdot 2^{j_t}, \mbox{ if } \epsilon_i=1\Big\}.$$
Note that
$${\rm Supp} \widehat{\Phi^{0}_{j_t, k}} \subset \Big\{|\xi_i|\leq \frac{4\pi}{3} \cdot 2^{j_t}, \forall i=1,\cdots, n\Big\}.$$
Hence
$$\left\langle B_{l,l',l'',\epsilon}(u,v), \Phi^{0}_{j_t,k}\right\rangle = 0.
$$

Then we consider the case where $\epsilon\neq 0$. We know
\begin{align*}
&a^{\epsilon'}_{j',k'}(t)= \left\langle B_{l,l',l'', \epsilon}(u,v), \Phi^{\epsilon'}_{j', k'}\right\rangle\\
&\;\;\;= \left\langle \int^{t}_{0} A^{t-s}_{l,l',l''} \Big\{\sum\limits_{\tiny \begin{array}{c} |j-j'|\leq 3, \epsilon, k, k''\\  j\geq 2+j_s\end{array}} u^{0}_{j-2,k}(s) v^{\epsilon}_{j,k''}(s) \Phi^{0}_{j-2,k}(x) \Phi^{\epsilon}_{j,k''}(x)\Big\} ds, \Phi^{\epsilon'}_{j', k'}\right\rangle.
\end{align*}
Hence
\begin{align*}
\left|a^{\epsilon'}_{j',k'}(t)\right| \lesssim &  2^{\frac{nj'}{2}+j'} \int ^{t}_{0} \sum\limits_{\tiny \begin{array}{c} |j-j'|\leq 3, \epsilon, k, k''\\  j\geq 2+j_s\end{array}}  \left|u^{0}_{j-2,k}(s)\right| \left|v^{\epsilon}_{j,k''}(s)\right|
e^{-c(t-s)2^{2j'} } ds \\
& \times (1+|4k-k''|)^{-N} (1+|2^{j-j'}k'-k''|)^{-N}.
\end{align*}
For all $j\geq 2+j_s$, we have
$\sup\limits_{k\in \mathbb{Z}^n }  |u^{0}_{j-2,k}(s)|
\lesssim  2^{j_s-\frac{nj}{2}}$ and
\begin{align*}
\left|a^{\epsilon'}_{j',k'}(t)\right|
\lesssim &  2^{j'} \int ^{t}_{0} 2^{ j_s} \sum\limits_{\tiny \begin{array}{c} |j-j'|\leq 3,\epsilon, k''\\  j\geq 2+j_s\end{array}}
\left|v^{\epsilon}_{j,k''}(s)\right| \left(1+|2^{j-j'}k'-k''|\right)^{-N}  \\
& \times e^{-c(t-s)2^{2j'} } ds .
\end{align*}
We decompose the integration for $t$ into two parts: (1) $\int ^{2^{-2(1+j_t)}}_{2^{-2-2j'} } $, (2) $\int ^{t}_{2^{-2(1+j_t)}}$.
We write the respective estimation by $a^{\epsilon',1}_{j',k'}(t)$ and $a^{\epsilon',2}_{j',k'}(t)$. Hence
$$\left|a^{\epsilon'}_{j',k'}(t)\right|\lesssim \left|a^{\epsilon',1}_{j',k'}(t)\right| + \left|a^{\epsilon',2}_{j',k'}(t)\right|.$$
For the first part, we have

\begin{align*}
&\left|a^{\epsilon',1}_{j',k'}(t)\right| \\
\lesssim &  2^{j'} e^{-ct2^{2j'} } \int ^{2^{-2(1+j_t)}}_{2^{-2-2j'} } 2^{ j_s} \sum\limits_{\tiny \begin{array}{c} |j-j'|\leq 3,\epsilon, k''\\  j\geq 2+j_s\end{array}}
\left|v^{\epsilon}_{j,k''}(s)\right| \left(1+|2^{j-j'}k'-k''|\right)^{-N}  ds\\
\lesssim &  2^{j'} e^{-ct2^{2j'} } \int ^{2^{-2(1+j_t)}}_{2^{-2-2j'} } 2^{ j_s} \left[\sum\limits_{\tiny \begin{array}{c} |j-j'|\leq 3\\  j\geq 2+j_s\end{array}}
\sum\limits_{\epsilon, k''}
\left|v^{\epsilon}_{j,k''}(s)\right|^{p} \left(1+\left|2^{j-j'}k'-k''\right|\right)^{-N}\right]^{\frac{1}{p}}  ds  .
\end{align*}

We choose $\tau$ satisfying $0<\tau<\min (m+\frac{1}{2}, \frac{1}{p'})$. By Cauchy-Schwartz inequality,
\begin{align*}
&|a^{\epsilon',1}_{j',k'}(t)|
\lesssim   2^{j'} e^{-ct2^{2j'} } \left(\int ^{2^{-2(1+j_t)}}_{2^{-2-2j'} } s^{-p'\tau} ds \right)^{\frac{1}{p'}}\\
&\;\;\;\times  \left[\int ^{2^{-2(1+j_t)}}_{2^{-2-2j'} } 2^{ pj_s}
  \sum\limits_{\tiny \begin{array}{c} |j-j'|\leq 3\\  j\geq 2+j_s\end{array}}
\sum\limits_{\epsilon, k''}
\left|v^{\epsilon}_{j,k''}(s)\right|^{p} \left(1+\left|2^{j-j'}k'-k''\right|\right)^{-N} s^{p\tau} ds\right]^{\frac{1}{p}}\\
&\;\;\;\;\;\;\;\;\;\;\lesssim  2^{j'} 2^{-2j'(\frac{1}{p'}-\tau)} e^{-ct2^{2j'} }\\
&\;\;\times\left[\int ^{2^{-2(1+j_t)}}_{2^{-2-2j'} } 2^{ pj_s} 2^{-2j_{s}p\tau}
\sum\limits_{\tiny \begin{array}{c} |j-j'|\leq 3\\  j\geq 2+j_s\end{array}}
\sum\limits_{\epsilon, k''}
\left|v^{\epsilon}_{j,k''}(s)\right|^{p} \left(1+|2^{j-j'}k'-k''|\right)^{-N}  ds\right]^{\frac{1}{p}} .
\end{align*}
Since $2m+1-2\tau>0$, we get
\begin{align*}
&\sum\limits_{j'\geq j_t}2^{2mp(j'-j_t)} 2^{j'p(\frac{n}{2}-1)}
\sum\limits_{\epsilon',k'} \left|a^{\epsilon',1}_{j',k'}(t)\right|^{p} \\
\lesssim & \sum\limits_{j'\geq j_t}2^{2mp(j'-j_t)} 2^{j'p\frac{n}{2}}
2^{-2pj'(\frac{1}{p'}-\tau)} e^{-cpt2^{2j'} }\\
&\times
\sum\limits_{\epsilon',k'}
\int ^{2^{-2(1+j_t)}}_{2^{-2-2j'} } 2^{ pj_s} 2^{-2j_{s}p\tau}
\sum\limits_{\tiny \begin{array}{c} |j-j'|\leq 3\\  j\geq 2+j_s\end{array}}
\sum\limits_{\epsilon, k''}
\left|v^{\epsilon}_{j,k''}(s)\right|^{p} \left(1+\left|2^{j-j'}k'-k''\right|\right)^{-N}  ds\\
\lesssim & \sum\limits_{j'\geq j_t}2^{2mp(j'-j_t)} 2^{j'p\frac{n}{2}}
2^{2pj'(\frac{1}{p'}-\tau)} e^{-cpt2^{2j'} }\\
&\times
\int ^{2^{-2(1+j_t)}}_{2^{-2-2j'} } 2^{ pj_s} 2^{-2j_{s}p\tau}
\sum\limits_{\tiny \begin{array}{c} |j-j'|\leq 3\\  j\geq 2+j_s\end{array}}
\sum\limits_{\epsilon, k''}
\left|v^{\epsilon}_{j,k''}(s)\right|^{p}  ds\\
\lesssim & \sum\limits_{j'\geq j_t}2^{2mp(j'-j_t)} 2^{j'p\frac{n}{2}}
2^{-2pj'(\frac{1}{p'}-\tau)} e^{-cpt2^{2j'} }\\
&\times
\sum\limits_{j_s=j_t} ^{1+j'}
2^{ pj_s} 2^{-2j_{s}p\tau}
2^{-2j_s} \sum\limits_{|j-j'|\leq 3} 2^{2mp(j_s-j)} 2^{-jp(\frac{n}{2}-1)}\\
\lesssim & \sum\limits_{j'\geq j_t}2^{2mp(j'-j_t)} 2^{j'p\frac{n}{2}}
2^{-2pj'(\frac{1}{p'}-\tau)} e^{-cpt2^{2j'} }
\sum\limits_{j_s=j_t} ^{1+j'}
2^{ pj_s} 2^{-2j_{s}p\tau}
2^{-2j_s}  2^{2mp(j_s-j')} 2^{-j'p(\frac{n}{2}-1)}.\\
\lesssim & \sum\limits_{j'\geq j_t}2^{2mp(j'-j_t)} 2^{j'p\frac{n}{2}}
2^{-2pj'(\frac{1}{p'}-\tau)} e^{-cpt2^{2j'} }
2^{ pj'} 2^{-2j'p\tau}
2^{-2j'}  2^{-j'p(\frac{n}{2}-1)} \\
\lesssim & \sum\limits_{j'\geq j_t}2^{2mp(j'-j_t)}  e^{-cpt2^{2j'} }
\lesssim
C.
\end{align*}


For the second part, we obtain
\begin{align*}
&\left|a^{\epsilon',2}_{j',k'}(t)\right|\\
 \lesssim & 2^{j'}  \int ^{t}_{2^{-2(1+j_t)}} 2^{ j_s}
\sum\limits_{\tiny \begin{array}{c} |j-j'|\leq 3,\epsilon, k''\\  j\geq 2+j_s\end{array}}
\left|v^{\epsilon}_{j,k''}(s)\right| \left(1+|2^{j-j'}k'-k''|\right)^{-N} e^{-c(t-s) 2^{2j'} }  ds\\
\lesssim & 2^{j'+j_t} \int ^{t}_{2^{-2(1+j_t)}}  \left[\sum\limits_{\tiny \begin{array}{c} |j-j'|\leq 3\\  j\geq 2+j_s\end{array}}
\sum\limits_{\epsilon, k''}
\left|v^{\epsilon}_{j,k''}(s)\right|^{p} \left(1+\left|2^{j-j'}k'-k''\right|\right)^{-N}\right]^{\frac{1}{p}}\\
&\times e^{-c(t-s) 2^{2j'} } ds \\
\lesssim & 2^{j'+j_t} \left[\int ^{t}_{2^{-2(1+j_t)}}   \sum\limits_{\tiny \begin{array}{c} |j-j'|\leq 3\\  j\geq 2+j_s\end{array}}
\sum\limits_{\epsilon, k''}
\left|v^{\epsilon}_{j,k''}(s)\right|^{p} \left(1+\left|2^{j-j'}k'-k''\right|\right)^{-N}
 e^{-c(t-s) 2^{2j'} } ds\right]^{\frac{1}{p}}\\
&\times  \left(\int ^{t}_{2^{-2(1+j_t)}} e^{-c(t-s) 2^{2j'} } ds \right)^{\frac{1}{p'}}\\
\lesssim & 2^{(1-\frac{2}{p'})j'+j_t} \\
&\,\,\,\times\left[\int ^{t}_{2^{-2(1+j_t)}}   \sum\limits_{\tiny \begin{array}{c} |j-j'|\leq 3\\  j\geq 2+j_s\end{array}}
\sum\limits_{\epsilon, k''}
|v^{\epsilon}_{j,k''}(s)|^{p} (1+|2^{j-j'}k'-k''|)^{-N}
e^{-c(t-s) 2^{2j'} } ds\right]^{\frac{1}{p}}.
\end{align*}
Hence
\begin{align*}
&\sum\limits_{j'\geq j_t}2^{2mp(j'-j_t)} 2^{j'p(\frac{n}{2}-1)}
\sum\limits_{\epsilon',k'} \left|a^{\epsilon',2}_{j',k'}(t)\right|^{p} \\
\lesssim & \sum\limits_{j'\geq j_t}2^{2mp(j'-j_t)} 2^{j'p (\frac{n}{2}-\frac{2}{p'})}
2^{pj_{t}}  \sum\limits_{\epsilon',k'} \int ^{t}_{2^{-2(1+j_t)}}   \sum\limits_{\tiny \begin{array}{c} |j-j'|\leq 3\\  j\geq 2+j_s\end{array}}
\sum\limits_{\epsilon, k''}
|v^{\epsilon}_{j,k''}(s)|^{p} \\
&\times (1+|2^{j-j'}k'-k''|)^{-N} e^{-c(t-s) 2^{2j'} } ds \\
\lesssim & \sum\limits_{j'\geq j_t}2^{2mp(j'-j_t)} 2^{j'p (\frac{n}{2}-\frac{2}{p'})}
2^{pj_{t}}   \int ^{t}_{2^{-2(1+j_t)}}   \sum\limits_{\tiny \begin{array}{c} |j-j'|\leq 3\\  j\geq 2+j_s\end{array}}
\sum\limits_{\epsilon, k''}
|v^{\epsilon}_{j,k''}(s)|^{p}  e^{-c(t-s) 2^{2j'} } ds
.
\end{align*}

We write the integration $\int ^{t}_{2^{-2(1+j_t)}}$ into two parts.
We take supremum of $\sum\limits_{\epsilon, k''}
|v^{\epsilon}_{j,k''}(s)|^{p}$ for $2^{-2(1+j_t)}\leq s< 2^{-2j_t}$ and $2^{-2j_t}\leq s <t$,
then we get the desired estimation.

\subsection{ Proof of (ii) of Lemma \ref{Le:4}}
We come to prove the following term belongs to $Y^{p}_{m}$: $$\int^{t}_{0} A^{t-s}_{l,l',l''} \left( \sum\limits_{j\geq 2+j_s} Q_j u Q_{j}v\right) ds.$$
We consider first the indices where $\epsilon=0$.
Note that
$$a^{0}_{j_t,k'}(t)= \int^{t}_{0}
\left\langle \sum\limits_{j\geq 2+j_s, k,k'',\epsilon,\epsilon'} u^{\epsilon}_{j,k}(s) v^{\epsilon'}_{j,k''}(s)
\Phi^{\epsilon}_{j,k}(x) \Phi^{\epsilon'}_{j,k''}(x), A^{t-s}_{l,l',l''} \Phi^{0}_{j_t, k'}\right\rangle ds.$$
By Lemma \ref{lem:es1}, we get
\begin{align*} &\left|a^{0}_{j_t,k'}(t)\right| \\
\lesssim & 2^{j_t+\frac{nj_t}{2}} \\
&\times\int^{t}_{0} \!\!\! \sum\limits_{j\geq 2+j_s, k,k'',\epsilon,\epsilon'}
\left|u^{\epsilon}_{j,k}(s)\right| \left|v^{\epsilon'}_{j,k''}(s)\right| \left(1+\left|k-k''\right|\right)^{-N} \left(1+|2^{j_t-j}k-k'|\right)^{-n-1}\!\!\! ds\\
\lesssim & 2^{j_t+\frac{nj_t}{2}} \sum\limits_{j'\geq 1+j_t } \int^{2^{-2j'}}_{2^{-2j'-2}}  \sum\limits_{j\geq 3+j'} a_{j,j',k'}(s)  ds
+ 2^{j_t+\frac{nj_t}{2}} \int^{t}_{2^{-2-2j_t}} \sum\limits_{j\geq 2+j_s} a_{j,k'}(s)  ds\\
\equiv & a^{0,1}_{j_t,k'}(t) + a^{0,2}_{j_t,k'}(t).
\end{align*}

Denote $\tilde{v}^{\epsilon'}_{j,k}(s)= \sum\limits_{k''}
|v^{\epsilon'}_{j,k''}(s)| (1+|k-k''|)^{-N}$.
This quantity satisfies the same estimation as $|v^{\epsilon'}_{j,k}(s)|$.
For $2^{-2j'-2}\leq s \leq 2^{-2j'}$, we deduce
\begin{align*}
a_{j,j',k'}(s)=& \sum\limits_{k,\epsilon,\epsilon'}
\left|u^{\epsilon}_{j,k}(s)\right| \left|\tilde{v}^{\epsilon'}_{j,k}(s)\right| \left(1+\left|2^{j_t-j}k-k'\right|\right)^{-n-1}\\
\lesssim & 2^{\frac{p-2}{p}(j-j_t)n} \left[\sum\limits_{k,\epsilon}
\left|u^{\epsilon}_{j,k}(s)\right| ^{p} \left(1+|2^{j_t-j}k-k'|\right)^{-n-1}\right]^{\frac{1}{p}}\\
&\times
\left[\sum\limits_{k,\epsilon'}
\left|\tilde{v}^{\epsilon'}_{j,k}(s)\right|^{p} \left(1+\left|2^{j_t-j}k-k'\right|\right)^{-n-1} \right]^{\frac{1}{p}}.
\end{align*}

Since $(\sum\limits_{k,\epsilon'}
|\tilde{v}^{\epsilon'}_{j,k}(s)|^{p})^{\frac{1}{p}}\lesssim 2^{2mj_s} 2^{(1-\frac{n}{2}-2m)j}$,
we know that
$$
a_{j,j',k'}(s)
\lesssim  2^{-\frac{p-2}{p}nj_t} 2^{2mj_s} 2^{(1+\frac{n}{2}-\frac{2n}{p}-2m)j}
\left[\sum\limits_{k,\epsilon}
\left|u^{\epsilon}_{j,k}(s)\right| ^{p} \left(1+\left|2^{j_t-j}k-k'\right|\right)^{-n-1}\right]^{\frac{1}{p}}.
$$
Take small positive $\delta$, we obtain
$$\sum\limits_{j\geq 3+j'} a_{j,j',k'}(s)\lesssim  \left(\sum\limits_{j\geq 3+j'} 2^{p\delta (j-j')} \left|a_{j,j',k'}(s)\right|^{p}\right)^{\frac{1}{p}}.$$

Choose small positive $\delta'$,
\begin{align*}
&2^{pj_t(\frac{n}{2}-1)} \sum\limits_{k'\in \mathbb{Z}^{n}} \left|a^{0,1}_{j_t,k'}(t)\right|^{p}\\
\lesssim & 2^{pnj_t} \sum\limits_{k'\in \mathbb{Z}^{n}} \sum\limits_{j'\geq 1+j_t }
2^{p\delta'(j'-j_t)} \left(\int^{2^{-2j'}}_{2^{-2j'-2}}  \sum\limits_{j\geq 3+j'} a_{j,j',k'}(s)  ds\right)^{p}\\
\lesssim & 2^{pnj_t} \sum\limits_{k'\in \mathbb{Z}^{n}} \sum\limits_{j'\geq 1+j_t } 2^{p\delta'(j'-j_t)}
2^{-2j'(p-1)} \int^{2^{-2j'}}_{2^{-2j'-2}}  \left(\sum\limits_{j\geq 3+j'} a_{j,j',k'}(s)\right)^{p}  ds       \\
\lesssim & 2^{pnj_t} \sum\limits_{k'\in \mathbb{Z}^{n}} \sum\limits_{j'\geq 1+j_t } 2^{p\delta'(j'-j_t)}
2^{-2j'(p-1)} \int^{2^{-2j'}}_{2^{-2j'-2}}  \sum\limits_{j\geq 3+j'} 2^{p\delta (j-j')} \left|a_{j,j',k'}(s)\right|^{p}  ds\\
\lesssim & 2^{pnj_t}  \sum\limits_{j'\geq 1+j_t } \sum\limits_{j\geq 3+j'} 2^{p\delta'(j'-j_t)} 2^{p\delta (j-j')}
2^{-2j'(p-1)} \int^{2^{-2j'}}_{2^{-2j'-2}}
\sum\limits_{k'\in \mathbb{Z}^{n}}\left|a_{j,j',k'}(s)\right|^{p}  ds.
\end{align*}

Further,
$A^{p,m}_{j,j'}= 2^{2mp(j-j')} 2^{jp(\frac{n}{2}-1)} \sup\limits_{1\leq t 2^{2j'}<4} \sum\limits_{(\epsilon,k)\in \Gamma}|a^{\epsilon}_{j,k}(t)|^{p} $.
For $2^{-2j'-2}\leq s< 2^{-2j'}$, we know
\begin{align*}
\sum\limits_{k'\in \mathbb{Z}^{n}}\left|a_{j,j',k'}(s)\right|^{p}
\lesssim &  2^{-(p-2) nj_t} 2^{2mpj'} 2^{(1+\frac{n}{2}-\frac{2n}{p}-2m)pj} \sum\limits_{k,\epsilon}
\left|u^{\epsilon}_{j,k}(s)\right| ^{p} \\ 
\lesssim &  2^{-(p-2) nj_t} 2^{4mpj'} 2^{(2-\frac{2n}{p}-4m)pj}
A^{p,m}_{j,j'+1}  .
\end{align*}
Thus,
\begin{align*}
&2^{pj_t(\frac{n}{2}-1)} \sum\limits_{k'\in \mathbb{Z}^{n}} \left|a^{0,1}_{j_t,k'}(t)\right|^{p}\\
\lesssim & 2^{2nj_t}  \sum\limits_{j'\geq 1+j_t } \sum\limits_{j\geq 3+j'} 2^{p\delta'(j'-j_t)} 2^{p\delta (j-j')}
2^{-2j'(p-1)} \int^{2^{-2j'}}_{2^{-2j'-2}}
2^{4mpj'} 2^{(2-\frac{2n}{p}-4m)pj}
A^{p,m}_{j,j'+1}   ds\\
\lesssim & 2^{2nj_t}  \sum\limits_{j'\geq 1+j_t } \sum\limits_{j\geq 3+j'} 2^{p\delta'(j'-j_t)} 2^{p\delta (j-j')}
2^{-2j'(p-1)} 2^{4mpj'} 2^{(2-\frac{2n}{p}-4m)pj} 2^{-2j'}
A^{p,m}_{j,j'+1}\\
= & 2^{2nj_t}  \sum\limits_{j'\geq 1+j_t } \sum\limits_{j\geq 3+j'} 2^{p\delta'(j'-j_t)-2nj'}
2^{(2+\delta-\frac{2n}{p}-4m)p(j-j')}
A^{p,m}_{j,j'+1}.
\end{align*}
We know
$\sum\limits_{j\geq j'} A^{p,m}_{j,j'}<\infty.$
For $0<\delta<\frac{2n}{p}+4m-2$ and $p\delta'<2n$,
we get
$$2^{pj_t(\frac{n}{2}-1)} \sum\limits_{k'\in \mathbb{Z}^{n}} \left|a^{0,1}_{j_t,k'}(t)\right|^{p} \lesssim C.
$$


For $2^{-2-2j_t}\leq s\leq t$, we obtain
\begin{align*}
&a_{j,k'}(s)= \sum\limits_{k,\epsilon,\epsilon'}
\left|u^{\epsilon}_{j,k}(s)\right| \left|\tilde{v}^{\epsilon'}_{j,k}(s)\right| \left(1+\left|2^{j_t-j}k-k'\right|\right)^{-n-1}\\
\lesssim & \left[\sum\limits_{k,\epsilon}
|u^{\epsilon}_{j,k}(s)| ^{p} (1+|2^{j_t-j}k-k'|)^{-n-1}\right]^{\frac{1}{p}}
\left[\sum\limits_{k,\epsilon'}
\left|\tilde{v}^{\epsilon'}_{j,k}(s)\right|^{p} \left(1+\left|2^{j_t-j}k-k'\right|\right)^{-n-1} \right]^{\frac{1}{p}}\\
&\times
\left[\sum\limits_{k}\left(1+\left|2^{j_t-j}k-k'\right|\right)^{-n-1}\right]^{\frac{p-2}{p}}\\
\lesssim & 2^{\frac{p-2}{p}(j-j_t)n} \left[\sum\limits_{k,\epsilon}
\left|u^{\epsilon}_{j,k}(s)\right| ^{p} \left(1+\left|2^{j_t-j}k-k'\right|\right)^{-n-1}\right]^{\frac{1}{p}}\\
&\times
\left[\sum\limits_{k,\epsilon'}
\left|\tilde{v}^{\epsilon'}_{j,k}(s)\right|^{p} \left(1+\left|2^{j_t-j}k-k'\right|\right)^{-n-1} \right]^{\frac{1}{p}}.
\end{align*}
Similarly, we have the desired estimation for $a^{0,2}_{j_t,k'}(t)$.


Then we consider the indices where $\epsilon'\neq 0$.
$$a^{\epsilon'}_{j',k'}(t)
= \int^{t}_{0} \left\langle \sum\limits_{j\geq 2+j_s, \epsilon,\tilde{\epsilon}, k,k''} u^{\epsilon}_{j,k} v^{\tilde{\epsilon}}_{j,k''} \Phi^{\epsilon}_{j,k}(x) \Phi^{\tilde{\epsilon}}_{j,k''}(x), A^{t-s,\beta}_{l,l',l''}\left \{\Phi^{\epsilon'}_{j', k'}\right\}\right\rangle ds.$$

Note that the following properties on the support of Fourier transform:
$$ {\rm Supp} \widehat{\Phi^{\epsilon}_{j,k} \Phi^{\tilde{\epsilon}}_{j,k''}} \subset
\Big\{|\xi_i|\leq \frac{16\pi}{3} \cdot 2^{j}, \forall i=1,\cdots,n\Big\},$$
$$ {\rm Supp} \widehat{\Phi^{\epsilon'}_{j', k'}} \subset
\Big\{|\xi_i|\leq \frac{4\pi}{3} \cdot 2^{j'}, \mbox{ if } \epsilon_i=0;
\frac{2\pi}{3} \cdot 2^{j'}\leq
|\xi_i|\leq \frac{8\pi}{3} \cdot 2^{j'}, \mbox{ if } \epsilon_i=1\Big\}.$$
Hence
\begin{align*}
|a^{\epsilon'}_{j',k'}(t)|
\lesssim & \int^{t}_{0}  \sum\limits_{j\geq \max( 2+j_s, j'-2) } \sum\limits_{\epsilon,\tilde{\epsilon}, k,k''} \left|u^{\epsilon}_{j,k}(s)\right| \left|v^{\tilde{\epsilon}}_{j,k''}(s)\right|
\left(1+|k-k''|\right)^{-N} \\
&\times \left(1+ \left|2^{j'-j}k-k'\right|\right)^{-N} 2^{(\frac{n}{2}+1)j'} e^{-c(t-s) 2^{2j'}}ds.
\end{align*}


Denote $\tilde{v}^{\tilde{\epsilon}}_{j,k}(s)= \sum\limits_{k''}
|v^{\tilde{\epsilon}}_{j,k''}(s)| (1+|k-k''|)^{-N}$. This quantity
 satisfies the same estimation as $|v^{\tilde{\epsilon}}_{j,k}(s)|$.
For $2^{-2j''}\leq s \leq 2^{2-2j''}$, since $$(\sum\limits_{k,\tilde{\epsilon}}
|\tilde{v}^{\tilde{\epsilon}}_{j,k}(s)|^{p})^{\frac{1}{p}}\lesssim 2^{2mj_s} 2^{(1-\frac{n}{2}-2m)j},$$
we deduce 
\begin{align*}
& a_{j,j',j'',k'}(s)\equiv \sum\limits_{k,\epsilon,\epsilon'}
\left|u^{\epsilon}_{j,k}(s)\right| \left|\tilde{v}^{\tilde{\epsilon}}_{j,k}(s)\right| \left(1+|2^{j'-j}k-k'|\right)^{-n-1}\\
&\lesssim  2^{\frac{p-2}{p}(j-j')n} \left[\sum\limits_{k,\epsilon}
|u^{\epsilon}_{j,k}(s)| ^{p} \left(1+\left|2^{j'-j}k-k'\right|\right)^{-n-1}\right]^{\frac{1}{p}}
\left[\sum\limits_{k,\tilde{\epsilon}}
|\tilde{v}^{\tilde{\epsilon}}_{j,k}(s)|^{p} (1+|2^{j'-j}k-k'|)^{-n-1} \right]^{\frac{1}{p}}\\
&\lesssim   2^{-\frac{p-2}{p}nj'} 2^{2mj_s} 2^{(1+\frac{n}{2}-\frac{2n}{p}-2m)j}
\left[\sum\limits_{k,\epsilon}
\left|u^{\epsilon}_{j,k}(s)\right| ^{p} \left(1+\left|2^{j'-j}k-k'\right|\right)^{-n-1}\right]^{\frac{1}{p}}.
\end{align*}


The proof for the case where $j_t\leq j'\leq 5+j_t$ is similar to the case where $j'\geq 6+j_t$.
We consider only the later. We divide the integration for $t$ into three parts:
$\sum\limits^{+\infty}_{j''=j'-3} \int^{2^{2-2j''}}_{2^{-2j''}},
\sum\limits^{j'-4}_{j''=2+j_t} \int^{2^{2-2j''}}_{2^{-2j''}},
\int^{t}_{2^{-(1+j_t)}}.$
Denote
\begin{align*}
\left|a^{\epsilon',1}_{j',k'}(t)\right|
&=  \sum\limits^{+\infty}_{j''=j'-3} \int^{2^{2-2j''}}_{2^{-2j''}}  \sum\limits_{j\geq \max( 2+j'', j'-2) }
a_{j,j',j'',k'}(s) 2^{(\frac{n}{2}+1)j'} e^{-c(t-s) 2^{2j'}}ds,\\
\left|a^{\epsilon',2}_{j',k'}(t)\right|
&=  \sum\limits^{j'-4}_{j''=2+j_t} \int^{2^{2-2j''}}_{2^{-2j''}}  \sum\limits_{j\geq \max( 2+j'', j'-2) }
a_{j,j',j'',k'}(s) 2^{(\frac{n}{2}+1)j'} e^{-c(t-s) 2^{2j'}}ds,\\
\left|a^{\epsilon',3}_{j',k'}(t)\right|
&=  \int^{t}_{2^{-(1+j_t)}}  \sum\limits_{j\geq \max( 2+j_s, j'-2) }
a_{j,j',j'',k'}(s) 2^{(\frac{n}{2}+1)j'} e^{-c(t-s) 2^{2j'}}ds.
\end{align*}
Hence
$$\left|a^{\epsilon'}_{j',k'}(t)\right|\lesssim \left|a^{\epsilon',1}_{j',k'}(t)\right| + \left|a^{\epsilon',2}_{j',k'}(t)\right| + \left|a^{\epsilon',3}_{j',k'}(t)\right|.$$

For the first, we obtain
\begin{align*}
& \left|a^{\epsilon',1}_{j',k'}(t)\right|
\lesssim  \sum\limits^{+\infty}_{j''=j'-3} \int^{2^{2-2j''}}_{2^{-2j''}}  \sum\limits_{j\geq  2+j'' }
a_{j,j',j'',k'}(s) 2^{(\frac{n}{2}+1)j'} e^{-ct 2^{2j'}}ds\\
\lesssim  & 2^{(\frac{n}{2}+1)j'} e^{-ct 2^{2j'}}  \sum\limits^{+\infty}_{j''=j'-3}
\int^{2^{2-2j''}}_{2^{-2j''}}  \left\{\sum\limits_{j\geq  2+j'' } 2^{p\delta (j-j'')}
\left|a_{j,j',j'',k'}(s) \right|^{p} \right\}^{\frac{1}{p}} ds\\
\lesssim  & 2^{(\frac{n}{2}+1)j'} e^{-ct 2^{2j'}}  \sum\limits^{+\infty}_{j''=j'-3} 2^{-\frac{2j''(p-1)}{p}}
\left\{\int^{2^{2-2j''}}_{2^{-2j''}}  \sum\limits_{j\geq  2+j'' } 2^{p\delta (j-j'')}
\left|a_{j,j',j'',k'}(s) \right|^{p}  ds \right\}^{\frac{1}{p}}\\
\lesssim  & 2^{(\frac{n}{2}+1)j'} e^{-ct 2^{2j'}} \\
&\times \left \{ \sum\limits^{+\infty}_{j''=j'-3} 2^{p\delta'(j''-j')} 2^{-2j''(p-1)}
\int^{2^{2-2j''}}_{2^{-2j''}}  \sum\limits_{j\geq  2+j'' } 2^{p\delta (j-j'')}
\left|a_{j,j',j'',k'}(s) \right|^{p}  ds \right\}^{\frac{1}{p}} .
\end{align*}

Accordingly,
\begin{align*}
&\sum\limits_{j'\geq j_t}2^{2mq(j'-j_t)} 2^{j'q(\frac{n}{2}-1)} \sum\limits_{(\epsilon',k')\in \Gamma}\left|a^{\epsilon',1}_{j',k'}(t)\right|^{p}\\
\lesssim & \sum\limits_{j'\geq j_t}2^{2mq(j'-j_t)} 2^{j'q(\frac{n}{2}-1)} 2^{(\frac{n}{2}+1)pj'} e^{-cpt 2^{2j'}}
\sum\limits_{(\epsilon',k')\in \Gamma}
\sum\limits^{+\infty}_{j''=j'-3} 2^{p\delta'(j''-j')} 2^{-2j''(p-1)}\\
&\times
\int^{2^{2-2j''}}_{2^{-2j''}}  \sum\limits_{j\geq  2+j'' } 2^{p\delta (j-j'')}
\left|a_{j,j',j'',k'}(s)\right |^{p}  ds\\
= & \sum\limits_{j'\geq j_t}2^{2mq(j'-j_t)} 2^{j'q(\frac{n}{2}-1)} 2^{(\frac{n}{2}+1)pj'} e^{-cpt 2^{2j'}}
\sum\limits^{+\infty}_{j''=j'-3} 2^{p\delta'(j''-j')} 2^{-2j''(p-1)}\\
&\times
\int^{2^{2-2j''}}_{2^{-2j''}}  \sum\limits_{j\geq  2+j'' } 2^{p\delta (j-j'')}
\sum\limits_{(\epsilon',k')\in \Gamma}
\left|a_{j,j',j'',k'}(s)\right |^{p}  ds\\
\lesssim & \sum\limits_{j'\geq j_t}2^{2mq(j'-j_t)} 2^{j'q(\frac{n}{2}-1)} 2^{(\frac{n}{2}+1)pj'} e^{-cpt 2^{2j'}}
\sum\limits^{+\infty}_{j''=j'-3} 2^{p\delta'(j''-j')} 2^{-2j''(p-1)}\\
&\times
\int^{2^{2-2j''}}_{2^{-2j''}}  \sum\limits_{j\geq  2+j'' } 2^{p\delta (j-j'')}
2^{-(p-2)nj_t} 2^{2mpj''} 2^{(1+\frac{n}{2}-\frac{2n}{p}-2m)pj}\\
&\times
\sum\limits_{(\epsilon',k')\in \Gamma} \sum\limits_{k,\epsilon}
\left|u^{\epsilon}_{j,k}(s)\right| ^{p} \left(1+|2^{j_t-j}k-k'|\right)^{-n-1} ds\\
\lesssim & \sum\limits_{j'\geq j_t}2^{2mq(j'-j_t)} 2^{j'q(\frac{n}{2}-1)} 2^{(\frac{n}{2}+1)pj'} e^{-cpt 2^{2j'}}
\sum\limits^{+\infty}_{j''=j'-3} 2^{p\delta'(j''-j')} 2^{-2j''(p-1)} \\
&\times
\int^{2^{2-2j''}}_{2^{-2j''}}  \sum\limits_{j\geq  2+j'' } 2^{p\delta (j-j'')}
2^{-(p-2)nj_t} 2^{2mpj''} 2^{(1+\frac{n}{2}-\frac{2n}{p}-2m)pj}
\sum\limits_{k,\epsilon} \left|u^{\epsilon}_{j,k}(s)\right| ^{p}  ds.
\end{align*}

Denote
$A^{p,m}_{j,j'}(u)= 2^{2mp(j-j')} 2^{jp(\frac{n}{2}-1)} \sup\limits_{1\leq t 2^{2j'}<4} \sum\limits_{(\epsilon,k)\in \Gamma}|u^{\epsilon}_{j,k}(t)|^{p} $.
For $2^{-2j''}\leq s< 2^{2-2j''},$ we know that
\begin{align*}
&\sum\limits_{j'\geq j_t}2^{2mq(j'-j_t)} 2^{j'q(\frac{n}{2}-1)} \sum\limits_{(\epsilon',k')\in \Gamma}\left|a^{\epsilon',1}_{j',k'}(t)\right|^{p}\\
\lesssim & \sum\limits_{j'\geq j_t}2^{2mq(j'-j_t)} 2^{j'q(\frac{n}{2}-1)} 2^{(\frac{n}{2}+1)pj'} e^{-cpt 2^{2j'}}
\sum\limits^{+\infty}_{j''=j'-3} 2^{p\delta'(j''-j')} 2^{-2j''(p-1)}
\int^{2^{2-2j''}}_{2^{-2j''}}   \\
&\times  \sum\limits_{j\geq  2+j'' } 2^{p\delta (j-j'')} 2^{-(p-2)nj_t} 2^{2mpj''} 2^{(1+\frac{n}{2}-\frac{2n}{p}-2m)pj}
2^{2mp(j''-j)} 2^{pj(1-\frac{n}{2})} A^{p,m}_{j,j''}(u)  ds\\
\lesssim & \sum\limits_{j'\geq j_t}2^{2mq(j'-j_t)} 2^{j'q(\frac{n}{2}-1)} 2^{(\frac{n}{2}+1)pj'} e^{-cpt 2^{2j'}}
\sum\limits^{+\infty}_{j''=j'-3} 2^{p\delta'(j''-j')} 2^{-2j''(p-1)}
2^{-2j''} \\
&\times  \sum\limits_{j\geq  2+j'' } 2^{p\delta (j-j'')}
2^{-(p-2)nj_t} 2^{2mpj''} 2^{(1+\frac{n}{2}-\frac{2n}{p}-2m)pj}
2^{2mp(j''-j)} 2^{pj(1-\frac{n}{2})} A^{p,m}_{j,j''}(u)\\
\lesssim & \sum\limits_{j'\geq j_t}2^{2mq(j'-j_t)} 2^{j'q(\frac{n}{2}-1)} 2^{(\frac{n}{2}+1)pj'} e^{-cpt 2^{2j'}}
\sum\limits^{+\infty}_{j''=j'-3} 2^{p\delta'(j''-j')} 2^{-2nj''} 2^{-(p-2)nj_t}\\
&\times  \sum\limits_{j\geq  2+j'' } 2^{p\delta (j-j'')} 2^{(4mp+2n-2p) (j''-j)}  A^{p,m}_{j,j''}(u).
\end{align*}

For $4mp+2n> p(2+\delta)$, we get
\begin{align*}
&\sum\limits_{j'\geq j_t}2^{2mq(j'-j_t)} 2^{j'q(\frac{n}{2}-1)} \sum\limits_{(\epsilon',k')\in \Gamma}\left|a^{\epsilon',1}_{j',k'}(t)\right|^{p}\\
\lesssim & \sum\limits_{j'\geq j_t}2^{2mq(j'-j_t)} 2^{j'q(\frac{n}{2}-1)} 2^{(\frac{n}{2}+1)pj'} e^{-cpt 2^{2j'}}
\sum\limits^{+\infty}_{j''=j'-3} 2^{p\delta'(j''-j')} 2^{-2nj''} 2^{-(p-2)nj_t}\\
\lesssim & \sum\limits_{j'\geq j_t}2^{2mq(j'-j_t)} 2^{j'q(\frac{n}{2}-1)} 2^{(\frac{n}{2}+1)pj'} e^{-cpt 2^{2j'}}
\sum\limits^{+\infty}_{j''=j'-3} 2^{p\delta'(j''-j')} 2^{-2nj''} 2^{-(p-2)nj_t}.
\end{align*}

For $p\delta'<2n$ and $p=q$, we deduce
\begin{align*}
&\sum\limits_{j'\geq j_t}2^{2mq(j'-j_t)} 2^{j'q(\frac{n}{2}-1)} \sum\limits_{(\epsilon',k')\in \Gamma}\left|a^{\epsilon',1}_{j',k'}(t)\right|^{p}\\
\lesssim & \sum\limits_{j'\geq j_t}2^{2mp(j'-j_t)} 2^{j'p(\frac{n}{2}-1)} 2^{(\frac{n}{2}+1)pj'} e^{-cpt 2^{2j'}} 2^{-2nj'} 2^{-(p-2)nj_t}\\
= & \sum\limits_{j'\geq j_t}2^{(2mp+pn-2n)(j'-j_t)}  e^{-cpt 2^{2j'}} \lesssim C.
\end{align*}

Similarly, we obtain

\begin{align*}
&\sum\limits_{j'\geq j_t}2^{2mp(j'-j_t)} 2^{j'p(\frac{n}{2}-1)} \sum\limits_{(\epsilon',k')\in \Gamma}\left|a^{\epsilon',2}_{j',k'}(t)\right|^{p}\\
\lesssim & \sum\limits_{j'\geq j_t}2^{2mp(j'-j_t)} 2^{j'p(\frac{n}{2}-1)} 2^{(\frac{n}{2}+1)pj'} e^{-cpt 2^{2j'}}
\sum\limits^{j'-4}_{j''=2+j_t} 2^{p\delta'(j''-j')} 2^{-2nj''} 2^{-(p-2)nj_t}\\
&\times  \sum\limits_{j\geq  2+j'' } 2^{p\delta (j-j'')} 2^{(4mp+2n-2p) (j''-j)}  A^{p,m}_{j,j''}\\
\lesssim & \sum\limits_{j'\geq j_t}2^{2mp(j'-j_t)} 2^{j'pn}  e^{-cpt 2^{2j'}}
\sum\limits^{j'-4}_{j''=2+j_t} 2^{p\delta'(j''-j')} 2^{-2nj''} 2^{-(p-2)nj_t}\\
\lesssim & \sum\limits_{j'\geq j_t}2^{(2m+n-\delta')p(j'-j_t)}   e^{-cpt 2^{2j'}}
\lesssim C.
\end{align*}

As for $a^{\epsilon',3}_{j',k'}(t)$, the proof is similar and more simple. We omit the details.

\subsection{ Proof of (iii) of Lemma \ref{Le:4}}
We estimate the norm of the third term $$ \int^{t}_{0} A^{t-s}_{l,l',l''} (  III_{u,v} ) ds .$$
We consider first the indices where $\epsilon=0$. Note that
$$
a^{0}_{j_t,k'}(t)
= \int^{t}_{0} \left\langle \sum\limits_{ k, k''} u^{0}_{j_s,k}(s) v^{0}_{j_s,k''} (s) \Phi^{0}_{j_s,k}(x) \Phi^{0}_{j_s,k''}(x),
A^{t-s}_{l,l',l''} ( \Phi^{0}_{j_t, k'})\right\rangle ds.
$$

Denote $|\tilde{v}^{0}_{j_s,k}(s) | =|v^{0}_{j_s,k''}(s) | (1+|k-k''|)^{-N}$.
Applying the estimates in Lemma \ref{kernel}, as we did  in  Lemma \ref{lem:es1}, we get

\begin{align*}
&|a^{0}_{j_t,k'}(t)| \\
\lesssim & 2^{j_t+ \frac{nj_t}{2}}\int^{t}_{0}  \sum\limits_{ k} \left|u^{0}_{j_s,k}(s)\right| \left|v^{0}_{j_s,k''}(s) \right| \left(1+|k-k''|\right)^{-N}
\left(1+\left|2^{j_t-j_s}k-k'\right|\right)^{-n-1} ds\\
\lesssim & 2^{j_t+ \frac{nj_t}{2}} \int^{t}_{0}
\left[\sum\limits_{k}
|u^{0}_{j_s,k}(s)| ^{p} (1+|2^{j_t-j_s}k-k'|)^{-n-1}\right]^{\frac{1}{p}}\\
&\times \left [\sum\limits_{k}
\left|\tilde{v}^{0}_{j_s,k}(s)\right|^{p} \left(1+\left|2^{j_t-j_s}k-k'\right|\right)^{-n-1} \right]^{\frac{1}{p}}
\left[\sum\limits_{k}\left(1+\left|2^{j_t-j_s}k-k'\right|\right)^{-n-1}\right]^{\frac{p-2}{p}} ds\\
\lesssim & 2^{j_t+ \frac{nj_t}{2}}  \int^{t}_{0}
\left[\sum\limits_{k} \left|u^{0}_{j_s,k}(s)\right| ^{p} \left(1+\left|2^{j_t-j_s}k-k'\right|\right)^{-n-1}\right]^{\frac{1}{p}} 2^{(1-\frac{n}{2})j_s} 2^{n(j_s-j_t)\frac{p-2}{p}} ds\\
=&  2^{j_t (1+\frac{2n}{p}- \frac{n}{2})}  \int^{t}_{0}
\left[\sum\limits_{k} \left|u^{0}_{j_s,k}(s)\right| ^{p} \left(1+\left|2^{j_t-j_s}k-k'\right|\right)^{-n-1}\right]^{\frac{1}{p}} 2^{(1+\frac{n}{2}-\frac{2n}{p})j_s}  ds.
\end{align*}

Take $p'= \frac{p}{p-1}$, $p'\tau<1, \delta= \frac{2}{p'}(1-p'\tau)\; and \; p\tau = p-1-\frac{p\delta}{2} $, we have

\begin{align*}
& \left|a^{0}_{j_t,k'}(t)\right| \\
\lesssim &   2^{j_t (1+\frac{2n}{p}- \frac{n}{2})}  \int^{t}_{0}
\left[\sum\limits_{k} \left|u^{0}_{j_s,k}(s)\right| ^{p} \left(1+\left|2^{j_t-j_s}k-k'\right|\right)^{-n-1}\right]^{\frac{1}{p}} 2^{(1+\frac{n}{2}-\frac{2n}{p})j_s}  ds\\
\lesssim & 2^{j_t (1+\frac{2n}{p}- \frac{n}{2}) } \left[\int^{t}_{0} 2^{( \frac{n}{2}+1-\frac{2n}{p})pj_s}
\sum\limits_{k} \left|u^{0}_{j_s,k}(s)\right| ^{p} \left(1+\left|2^{j_t-j_s}k-k'\right|\right)^{-n-1} s^{p\tau} ds \right]^{\frac{1}{p}} \\
&\times  \left(\int^{t}_{0} s^{-p'\tau} ds\right)^{\frac{1}{p'}} \\
\lesssim & 2^{j_t (1+\frac{2n}{p}- \frac{n}{2}-\delta ) }\\
&\times \left[\int^{t}_{0} 2^{( \frac{n}{2}+1-\frac{2n}{p})pj_s - 2(p-1) j_s + p\delta j_s }
\sum\limits_{k} \left|u^{0}_{j_s,k}(s)\right| ^{p} \left(1+\left|2^{j_t-j_s}k-k'\right|\right)^{-n-1}  ds \right]^{\frac{1}{p}}\\
= &  2^{j_t (1+\frac{2n}{p}- \frac{n}{2}-\delta ) } \\
&\times\left[\int^{t}_{0} 2^{( \frac{n}{2}-1-\frac{2n}{p} + \delta)pj_s +2 j_s  }
\sum\limits_{k} \left|u^{0}_{j_s,k}(s)\right| ^{p} \left(1+\left|2^{j_t-j_s}k-k'\right|\right)^{-n-1}  ds\right ]^{\frac{1}{p}} .
\end{align*}
Hence
\begin{align*}
&2^{pj_t(\frac{n}{2}-1)} \sum\limits_{k'\in \mathbb{Z}^{n}} \left|a^{0}_{j_t,k'}(t)\right|^{p}\\
\lesssim & 2^{pj_t(\frac{n}{2}-1)} \sum\limits_{k'\in \mathbb{Z}^{n}}
2^{pj_t (1+\frac{2n}{p}- \frac{n}{2}-\delta) }\\
&\times \int^{t}_{0} 2^{(\delta+ \frac{n}{2}-1-\frac{2n}{p})pj_s + 2j_s}
\sum\limits_{k} \left|u^{0}_{j_s,k}(s)\right| ^{p} \left(1+\left|2^{j_t-j_s}k-k'\right|\right)^{-n-1} ds\\
\lesssim & 2^{pj_t (\frac{2n}{p}- \delta) } \int^{t}_{0} 2^{(\delta+ \frac{n}{2}-1-\frac{2n}{p})pj_s + 2j_s}
\sum\limits_{k} \left|u^{0}_{j_s,k}(s)\right| ^{p}  ds\\
\lesssim & 2^{pj_t (\frac{2n}{p}- \delta) } \sum\limits_{j_s\geq j_t}  2^{(\delta+ \frac{n}{2}-1-\frac{2n}{p})pj_s}
2^{pj_s (1-\frac{n}{2})} \lesssim C.
\end{align*}

Then we consider the indices where $\epsilon\neq 0$.
$$a^{\epsilon}_{j',k'}(t)= \int^{t}_{0} \left\langle \sum\limits_{ k, k''} u^{0}_{j_s,k}(s) v^{0}_{j_s,k''}(s) \Phi^{0}_{j_s,k}(x) \Phi^{0}_{j_s,k''}(x),
A^{t-s}_{l,l',l''}( \Phi^{\epsilon}_{j', k'}) \right\rangle ds.$$

Note that ${\rm Supp} \widehat{\Phi^{0}_{j_s,k}} \subset \Big\{|\xi_i|\leq \frac{4\pi}{3} \cdot 2^{j_s}, \forall i=1,\cdots, n\Big\}$ implies
$${\rm Supp} \widehat{III_{u,v}} \subset \Big\{|\xi_i|\leq \frac{8\pi}{3} \cdot 2^{j_s}, \forall i=1,\cdots, n\Big\}.$$
Further,
$$ {\rm Supp} \widehat{\Phi^{\epsilon'}_{j', k'}} \subset
\Big\{|\xi_i|\leq \frac{4\pi}{3} \cdot 2^{j'}, \mbox{ if } \epsilon_i=0;
\frac{2\pi}{3} \cdot 2^{j'}\leq
|\xi_i|\leq \frac{8\pi}{3} \cdot 2^{j'}, \mbox{ if } \epsilon_i=1\Big\}.$$
Only the case $j_s \geq j'-1$,
$$\left\langle  \Phi^{0}_{j_s,k}(x) \Phi^{0}_{j_s,k''}(x),
A^{t-s}_{l,l',l''} ( \Phi^{\epsilon}_{j', k'}) \right\rangle \neq 0.$$

Depending on the value of $j'$, we estimate it in two cases:
(i) $j_t\leq j'\leq 2+j_t$; (ii) $j'\geq 3+j_t$.
By similarity, we consider only the later.
For $j'\geq 3+j_t$,
\begin{align*}
&\left|a^{\epsilon}_{j',k'}(t)\right|\lesssim  2^{j'+ \frac{nj'}{2}} \int^{2^{-2j'}}_{0} \sum\limits_{ k} \left|u^{0}_{j_s,k}(s) \right| \left|v^{0}_{j_s,k''}(s)\right|(1+|k-k''|)^{-N}
\\
&\;\;\times \left(1+\left|2^{j'-j_s}k-k'\right|\right)^{-N}  e^{-ct 2^{2j'}} ds\\
\lesssim &   2^{j'+ \frac{nj'}{2}} e^{-ct 2^{2j'}} \int^{2^{-2j'}}_{0}  \sum\limits_{ k} \left|u_{j_s,k} \right| \left|\tilde{v}_{j_s,k}\right|
\left(1+\left|2^{j'-j_s}k-k'\right|\right)^{-N} ds\\
\lesssim &   2^{j'+ \frac{nj'}{2}} e^{-ct 2^{2j'}} \int^{2^{-2j'}}_{0}
\left[\sum\limits_{k} \left|u^{0}_{j_s,k}(s)\right| ^{p} \left(1+\left|2^{j'-j_s}k-k'\right|\right)^{-n-1}\right]^{\frac{1}{p}}\\
&\times  \left[\sum\limits_{k}
\left|\tilde{v}^{0}_{j_s,k}(s)\right|^{p} \left(1+\left|2^{j'-j_s}k-k'\right|\right)^{-n-1} \right]^{\frac{1}{p}}
\left[\sum\limits_{k}\left(1+|2^{j'-j_s}k-k'|\right)^{-n-1}\right]^{\frac{p-2}{p}} ds \\
\lesssim & 2^{j' (1+\frac{2n}{p}- \frac{n}{2})} e^{-ct 2^{2j'}} \int^{2^{-2j'}}_{0}
\left[\sum\limits_{k} \left|u^{0}_{j_s,k}(s)\right| ^{p} \left(1+\left|2^{j'-j_s}k-k'\right|\right)^{-n-1}\right]^{\frac{1}{p}} 2^{(1+\frac{n}{2}-\frac{2n}{p})j_s}  ds\\
\lesssim & 2^{j' (1+\frac{2n}{p}- \frac{n}{2}-\delta)} e^{-ct 2^{2j'}} \\
&\times\left[\int^{2^{-2j'}}_{0}
2^{( \frac{n}{2}-1-\frac{2n}{p} + \delta)pj_s +2 j_s  }
\sum\limits_{k} \left|u^{0}_{j_s,k}(s)\right| ^{p} \left(1+\left|2^{j_t-j_s}k-k'\right|\right)^{-n-1}  ds \right]^{\frac{1}{p}}
. \end{align*}
That is to say,
\begin{align*}
&\sum\limits_{j'\geq 3+ j_t}2^{2mp(j'-j_t)} 2^{j'p(\frac{n}{2}-1)} \sum\limits_{(\epsilon',k')\in \Gamma}\left|a^{\epsilon'}_{j',k'}(t)\right|^{p}\\
\lesssim & \sum\limits_{j'\geq j_t}2^{2mp(j'-j_t)} 2^{j'p(\frac{n}{2}-1)} 2^{pj' (1+\frac{2n}{p}- \frac{n}{2}-\delta)} e^{-cpt 2^{2j'}}\\
&\times \sum\limits_{(\epsilon',k')\in \Gamma}
\int^{2^{-2j'}}_{0}
2^{( \frac{n}{2}-1-\frac{2n}{p} + \delta)pj_s +2 j_s  }
\sum\limits_{k} \left|u^{0}_{j_s,k}(s)\right| ^{p} \left(1+\left|2^{j_t-j_s}k-k'\right|\right)^{-n-1}  ds\\
\lesssim & \sum\limits_{j'\geq j_t}2^{2mp(j'-j_t)}  2^{pj' (\frac{2n}{p}-\delta)} e^{-cpt 2^{2j'}}\int^{2^{-2j'}}_{0}
2^{( \frac{n}{2}-1-\frac{2n}{p} + \delta)pj_s +2 j_s  }
\sum\limits_{k} \left|u^{0}_{j_s,k}(s)\right| ^{p}   ds\\
\lesssim & \sum\limits_{j'\geq j_t}2^{2mp(j'-j_t)}   e^{-cpt 2^{2j'}} \lesssim C.
\end{align*}

\section{Single norm and global smooth solution}\label{sec6}

For initial data $a(x)$ belong to some initial space $\big(X(\mathbb{R}^n)\big)^{n}$,
usually the mild solution of (\ref{eqn:ns}) is found in some solution space $\big(Y(\mathbb{R}^{1+n}_{+})\big)^n$
near the solution $e^{t\Delta} a(x)$ of the heat equation
which is obtained by using the following method via the iteration process:
$$
\begin{cases}
u^{(0)}(t,x)=e^{t\Delta} a(x);\\
u^{(j+1)}(t,x)=u^{(0)}(t,x)-B(u^{(j)},u^{(j)})(t,x)\quad\hbox{for}\quad
j=0,1,2,....
\end{cases}
$$
The iteration process has  the following properties:

(i) If $a$ is in $\big(X(\mathbb{R}^n)\big)^{n}$, then
$e^{t\Delta}a(x)$ belongs to $\big(Y(\mathbb{R}^{1+n}_{+})\big)^{n}$.

(ii) The bilinear operator
$$
\begin{array}{rl}
B(u,v)=&\int^{t}_{0} e^{(t-s)\Delta} \mathbb{P}\nabla
(u\otimes v) ds\end{array}
$$
is bounded from $(Y(\mathbb{R}^{1+n}_{+}))^n \times
(Y(\mathbb{R}^{1+n}_{+}))^n $ to $(Y(\mathbb{R}^{1+n}_{+}))^n $.

With the initial data being small in $\big(X(\mathbb{R}^n)\big)^{n}$,
we get a contraction mapping on  $\big (Y(\mathbb{R}^{1+n}_{+})\big)^{n}$.
The fixed point theorem implies that
there exists a unique mild solution of \eqref{eqn:ns} in
$\big(Y(\mathbb{R}^{1+n}_{+})\big)^{n}$.
Regarding the well-posedness results of the Navier-Stokes equations,
the differences in the way we deal with it are the choice of initial value space,
the choice of solution space reflecting the structure of the solution,
and the method of proving the continuity of (i) and (ii) above.

\begin{remark}
The main purpose of this paper is to explore the internal harmonic analysis structure of the solution of the equations \eqref{eqn:ns}.
The initial data space $\big(X(\mathbb{R}^{n})\big)^n$
is taken as the usual critical Besov spaces $\dot{B}^{\frac{n}{p}-1,p}_{p}$.
But the structure of the solution space is different.
Usually, we find solution in the subspace of $L^{\infty}((\dot{B}^{\frac{n}{p}-1,p}_{p})^{n})$.
Further, all the previous solution spaces belong to the intersection of two spaces with different norms.
The solution space here is the space defined by
the decaying property of the norm of the derivative of a given order over time.
With the help of {\bf Gevrey class} in \cite{B},
our solution can also be given in terms of analytic or decaying properties of the infinitely smooth norm.
\end{remark}

Based on the present results, there are at least two differences between the fractional equation and the classical equation.
For the latter, firstly, it is still not known whether Bloch space is well-posed ;
Secondly, the solution space is the intersection of two different norm spaces.
At first we did not expect that
the decaying property of a single norm with respect to time
based on the space-time structure would match the classical equation.
Building on the new ideas in the previous sections,
our other ideas about the well-posedness Theorem \ref{mthmain} are also based on the traditional idea of mild solutions.
Let us prove  Theorem \ref{mthmain}.

\begin{proof} (i) Given  $n<p<\infty$ and $m>0 $. By Theorem \ref{th:B-to-Y}, we know, if
$f\in \dot{B}^{\frac{n}{p}-1,p}_{p}$, then $ e^{t\Delta} f \in Y^{p}_{m}.$

(ii) For all $l, l',l''=1,\cdots,n$, let
$B_{l}(u,v)$ and $B_{l,l',l''}(u,v)$ be defined in the equations (\ref{41eq}) and (\ref{42eq}) respectively.
By Theorem \ref{th:4.1}, we know that all $B_{l}(u,v)$,
$B_{l,l',l''}(u,v)$ are bounded from $Y^{p}_{m}\times Y^{p}_{m}$ to $Y^{p}_{m}$.
These imply that the bilinear operator
\begin{equation*}
B(u,v)= \int^{t}_{0} e^{(t-s)\Delta}
\mathbb{P}\nabla (u\otimes v) ds
\end{equation*}
is bounded from $(Y^{p}_{m})^{n}\times (Y^{p}_{m})^{n}$
to $(Y^{p}_{m})^{n}$.

(iii) Hence, by Picard's contraction principle, for small initial data in $(\dot{B}^{\frac{n}{p}-1,p}_{p})^{n}$,
there exists a unique solution in $(Y^{p}_{m})^{n}$.
\end{proof}

{\bf Acknowledgements.}
The first author would like to thank Professor Huijiang Zhao for his invaluable discussions and suggestions.


\end{document}